\documentclass[a4paper,10pt]{amsart}
\usepackage[english]{certus}
\usepackage{a4wide}

\title[Sofic boundaries and coarse geometry of sofic approximations]{Sofic boundaries of groups and coarse geometry of sofic approximations}
\author{Vadim Alekseev}
\address{Vadim Alekseev, Technische Universit\"{a}t Dresden, Fachrichtung Mathematik, Institut f\"{u}r Geometrie, 01062, Dresden, Deutschland}
\email{vadim.alekseev@tu-dresden.de}

\author{Martin Finn--Sell}
\address{Martin Finn--Sell, Universit{\"a}t Wien, Fakult\"{a}t f\"{u}r Mathematik, Oskar-Morgenstern-Platz 1,   1090 Wien, \"{O}sterreich }
\email{martin.finn-sell@univie.ac.at}


\begin{document}
\onehalfspace

\begin{abstract}
Sofic groups generalise both residually finite and amenable groups, and the concept is central to many important results and conjectures in measured group theory. We introduce a topological notion of a sofic boundary attached to a given sofic approximation of a finitely generated group and use it to prove that coarse properties of the approximation (property A, asymptotic coarse embeddability into Hilbert space, geometric property (T)) imply corresponding analytic properties of the group (amenability, a-T-menability and property (T)), thus generalising ideas and results present in the literature for residually finite groups and their box spaces. Moreover, we generalise coarse rigidity results for box spaces due to Kajal Das, proving that coarsely equivalent sofic approximations of two groups give rise to a uniform measure equivalence between those groups. Along the way, we bring to light a coarse geometric view point on ultralimits of a sequence of finite graphs first exposed by J\'{a}n \v{S}pakula and Rufus Willett, as well as proving some bridging results concerning measure structures on topological groupoid Morita equivalences that will be of interest to groupoid specialists.
\end{abstract}

\maketitle

\section{Introduction}

Finite approximation of infinite objects is a fundamental tool in the modern mathematician's toolkit, and it has been used to great effect in the authors' favourite areas of mathematics: in the realm of operator algebras the notions of nuclearity, exactness and quasidiagonality for \Cs-algebras \cite{MR3418247,TikuisisQuasidiagonalityofnuclearCalgebras2015,MR1437044}, and the corresponding notion of hyperfiniteness for von Neumann algebras \cite{MurrayOn1943} have given rise to the classification programs of $C^{*}$-algebras \cite{ElliottOn1976,kirchberg-classification-1999} and von Neumann algebras \cite{ConnesClassification1976}. Their natural group theoretic counterpart is amenability.

The aforementioned types of approximation are quite strong and therefore restrictive: they correspond to the ``amenable world'' of groups and operator algebras. While interesting and beautiful in its own right, it does not encompass many natural and important examples in group theory and operator algebras -- say, the free groups and operator algebraic objects related to them. However, one would like to extend the idea of finitary approximation as well beyond amenability. In the realm of operator algebras, such an approximation was suggested by Alain Connes in \cite{ConnesClassification1976} and lead to the famous \emph{Connes Embedding Conjecture}. By the remarkable work of Eberhard Kirchberg \cite{KirchbergOn1993} it was shown to be equivalent to the so-called QWEP conjecture for \Cs-algebras.

What one sees by studying the above is a relaxation of \textit{algebra homomorphisms} to maps that are \textit{approximately} homomorphisms. This suggests a more general notion of finite approximation should exist for groups when we allow for a \textit{metric} on the finite set on which we attempt to approximate. This leads to the definition of a \textit{sofic} group.

To make sense of what an ``approximate" map to a finite group is, one chooses finite symmetric groups as targets and equips them with the normalised Hamming distance. A group $\Gamma$ is sofic if it is possible to find approximations of arbitrary finite subsets of $\Gamma$ in symmetric groups $\Sym(X)$ that are approximately injective and approximately multiplicative with respect to this distance. A countable collection $\mc X$ of finite sets $X_{i}$ that witness stronger and stronger approximations for an exhaustion of the group $\Gamma$ is a \textit{sofic approximation} of $\Gamma$. Examples of sofic groups include amenable groups and residually finite discrete groups. Sofic groups were introduced by Mikhail Gromov \cite{MR1694588} in his work on Gottschalk's surjunctivity conjecture, and expanded on (and named by) Benjamin Weiss in \cite{MR1803462}. Since then they have played a fundamental role in research in dynamical systems.

% , notably being the largest class of groups for which there are notion of entropy is well defined \cite{MR2552252} and the notion of mean dimension has been extended \cite{MR3077882}. These ideas have since been extended to discrete measured groupoids \cite{MR3130315}.

The purpose of this paper is to introduce a general technique for studying sofic approximations of groups from the coarse geometric point of view and to give a mechanism for transferring topological (in this context, coarse geometric) properties from the approximation back to the group. The vessel we use to complete this journey is coarse geometric in nature and was initially introduced by George Skandalis, Jean-Louis Tu and Guoliang Yu in \cite{MR1905840}, where a \textit{topological groupoid} was constructed to emulate the role of a group in certain aspects of the Baum--Connes conjecture for metric spaces. The second author of this paper studied this groupoid and certain of its \textit{reductions} in \cite{MR3197659} and \cite{MR3266245} in the context of \textit{box spaces} associated to residually finite discrete groups. 

A \textit{box space} associated to a residually finite discrete group $\Gamma$ and a chain of subgroups $\lbrace N_{i} \rbrace_{i}$ is a metric space, denoted $\square \Gamma$, constructed from the Cayley graphs of the finite quotients $\Gamma/N_{i}$. % ; this provides a basic example of a \textit{space of graphs},  Section \ref{sect:coarse}.
 This is a particular example of a sofic approximation of a residually finite group.
 
 Box spaces can be a powerful tool, both to differentiate between coarse properties (as in \cite{MR2899681}) and to provide a finite dimensional test for analytic properties of the group $\Gamma$. Notably, the following correspondences between coarse geometric properties of the box space and analytic properties of the group are known:
\begin{itemize}
\item $\square \Gamma$ has Property A if and only if  $\Gamma$ is amenable \cite[Proposition 11.39]{MR2007488};
\item $\square \Gamma$ has an asymptotic coarse embedding (or a fibred coarse embedding)
%\footnote{Originally this equivalence would have been phrased using a \textit{fibred coarse embedding into Hilbert space} (see \cite{MR3116568,MR3105001}). However, Rufus Willett \cite{MR3346926} introduced the more natural \textit{asymptotic coarse embedding}, which in this instance is equivalent to a fibred coarse embedding and easier to formulate.}
 into Hilbert space if and only if $\Gamma$ is a-T-menable \cite{MR3346926,MR3266245,MR3116568,MR3105001};
\item $\square \Gamma$ has geometric property (T) if and only if $\Gamma$ has property (T) \cite{MR3246936}.
\end{itemize}
% In many instances, these results were broadened by considering the space of marked groups in \cite{Mimura-Sako,MR3342685}. 

The method presented in \cite{MR3266245} for producing these results was to associate to any given box space $\square \Gamma$ a topological boundary that admits a free $\Gamma$-action -- this boundary action is a particular component of the coarse groupoid of Skandalis--Tu--Yu. The main idea in this paper is to generalise this procedure to a sofic approximation of a sofic group, but in this setting the counting measures on each ``box" will play a fundamental role. More precisely, we associate to a given sofic approximation a topological groupoid that we call the \emph{sofic coarse boundary groupoid}. The base space of this groupoid -- the \emph{sofic boundary} -- is constructed from the ``box space'' of graphs coming from the sofic approximation. It carries a natural invariant measure coming from the counting measure on the graphs and has a nice closed saturated subset $Z$ of full measure -- the \emph{core} of the sofic boundary -- restricted to which, the sofic coarse boundary groupoid turns out to be  a crossed product by an action of $\Gamma$ as in the traditional box space case. This allows us to prove:

\begin{Thm}\label{thm:intro}
Let $\Gamma$ be a sofic group, $\mc X$ a sofic approximation of $\Gamma$, and $X$ be the space of graphs constructed from $\mc X$. Then:
\begin{enumerate}
\item If $X$ has property A then  $\Gamma$ is amenable (Theorem \ref{Thm:amenable});
\item If $X$ admits an asymptotic coarse embedding into Hilbert space, then $\Gamma$ is a-T-menable (Theorem \ref{Thm:a-t-men});
\item If $X$ has boundary geometric property (T) then $G$ has property (T) (Theorem \ref{Thm:T}).
\end{enumerate}
\end{Thm}
% These results are explained in the corresponding part of Section \ref{sect:main-theorems} associated to each of these properties (and are more complicated, as they concern representation theory of groupoids -- the basics of which are included in Appendix \ref{app:rep}). In Section \ref{sect:LEF},

At this point, it is natural to ask about the converse statements. There appears to be little hope of establishing them in full generality, the main technical reason being that the core of the sofic boundary is a proper subset of it, and there is no control of what happens on the complement. We explain this issue in more detail in the final section of the paper.

However, if the sofic approximation is coming from the group being \textit{locally embeddable into a finite group} (or briefly an LEF group), the core is the entire boundary, which allows us to recover the converse to the above statements, thus reproving the known results about LEF groups from the literature \cite{Mimura-Sako,MR3342685}. 
% We give two examples of groups that are LEF, which have somewhat interesting approximations as a consequence of these results. The first is the full group of a Cantor minimal system \cite{MR3071509}, which are LEF \cite{MR3241829}, and have natural infinite amenable subgroups. What makes them interesting is they are simple, and so a traditional ``box space'' is not useful here. The second example is the a-T-menable, non-$C^{*}$-exact group of Osajda \cite{Osajda-nonexact}, which is also LEF (as it is a direct limit of CAT(0)-cubical hyperbolic groups).

Transitioning from coarse invariants (that are \textit{topological} invariants of a groupoid) to measurable invariants, we begin to investigate the question: to what extent a sofic approximation is a ``coarse invariant'' of the sofic group? To this end, we were able to prove the following:

\begin{Thm}(Theorem \ref{thm:coarse-equiv})
Let $\Gamma$, $\Lambda$ be sofic groups with sofic approximations $\mc X$ and $\mc Y$ respectively. Let $X_{\mc X}$ and $X_{\mc Y}$ be their associated spaces of graphs. If $X_{\mc X}$ and $X_{\mc Y}$ are coarsely equivalent, then $\Gamma$ and $\Lambda$ are quasi-isometric and uniformly measure equivalent.
\end{Thm}

This theorem generalises part of the work in \cite{Khukhro-Valette}, and the main result of \cite{Das-box} to the case that $\Gamma$ and $\Lambda$ are sofic, as opposed to residually finite, and the technique is completely different -- we construct a \textit{Morita equivalence} bispace for the sofic coarse boundary groupoids. This bispace looks very much like the topological coupling introduced by Gromov in his dynamic classification of quasi-isometries between groups. Given appropriate measures on the groupoids, we construct a measure on the bispace, which turns the topological Morita equivalence into a measurable one -- and this allows us to deduce the uniform measure equivalence combining the topological and measure-theoretic properties of sofic coarse boundary groupoids. As was pointed out in \cite{Das-box}, by combining a result of Damien Gaboriau \cite[Theorem 6.3]{MR1953191} with Theorem \ref{thm:coarse-equiv} we are able to conclude facts concerning the rigidity of $\ell^{2}$-Betti numbers of sofic groups with coarsely equivalent approximations:

\begin{Cor}
If $\Gamma$ and $\Lambda$ are finitely generated sofic groups with coarsely equivalent sofic approximations, then their $\ell^{2}$-Betti numbers are proportional.
\end{Cor}

The downside of the topological groupoid we construct to settle the above questions is that the unit space is not second countable, therefore not metrizable (and thus not a \textit{standard} as a probability space). We remedy this situation by providing a recipe for constructing many different second countable versions of the groupoid using ideas from \cite{MR1905840,MR2419901}. The following should be considered as a topological result in line with the standartisation theorem for measurable actions proved by Alessandro Carderi in \cite[Theorem A]{carderi-ultraproducts}.

\begin{Thm}
Let $\Gamma$ be a sofic group, $\mc X$ a sofic approximation of $\Gamma$, $X$ the associated total space of the family of graphs attached to $\mc X$ and $Z\subset X$ the core of a sofic approximation. Then there exists a second countable \'etale, locally compact, Hausdorff topological groupoid $\mc G$ with following properties:
\begin{enumerate}
\item the base space $\mc G^{(0)}\eqqcolon \widehat X$ is a compactification of $X$ (in particular, it's a quotient of $\beta X$ through a quotient map $p\colon \beta X \to \widehat X$),
\item $p(Z)\subset \partial \widehat X$ is invariant and satisfies $\mc G|_{p(Z)} \cong p(Z)\rtimes \Gamma$. As a consequence, we have an almost everywhere isomorphism
\[
(\mc G|_{\partial \widehat X}, \nu_{p_{*}\mu}) \rightarrow (\widehat X,p_{*}\mu)\rtimes \Gamma.
\]
\end{enumerate}
\end{Thm}

As an example of this process, we construct the minimal topological groupoid introduced in \cite{MR2966663} for a residually finite discrete group and a corresponding Farber chain of finite index subgroups. 

The paper is organised as follows. In Section \ref{sect:prelim} we recapitulate the necessary definitions and results both from the theory of sofic group approximations and groupoids arising from coarse geometry. Section \ref{sect:sofic-coarse-bd-groupoid} introduces our main player, the sofic coarse boundary groupoid associated with a fixed sofic approximation of a group and studies its properties; in particular, we introduce the core of a sofic approximation as the closure of the ``good set'' in the approximating graphs. Section \ref{sect:main-theorems} is devoted to the proof of the main Theorem \ref{thm:intro} and its converse in the case of an LEF group. Finally, in Section \ref{sect:coarseequiv} we prove that coarse equivalence of two sofic approximations implies quasi-isometry and uniform measure equivalence of groups (Theorem \ref{thm:coarse-equiv}). In the last section we discuss some related open questions that might be of interest for further investigation.

\section{Preliminaries}\label{sect:prelim}
In this section we introduce the necessary definitions, facts and references for coarse groupoids and sofic groups.

\subsection{Groupoids from coarse geometry}\label{sect:groupoids}
We recapitulate some particular examples of groupoids that appear later in the paper. For a basic introduction to \'etale groupoids we recommend \cite{MR2419901}, for their representation theory \cite{MR2966043} and finite approximation properties \cite{MR1799683}. We also suggest the collected references of \cite{MR1905840}, \cite{MR2007488} and \cite{Spakula-Willett} for the notion of coarse groupoid and its properties. 

\begin{Ex}\label{Ex:TransGrp}
Let $X$ be a topological $\Gamma$-space. Then the \textit{transformation groupoid} associated to this action is given by the data $X \times G \rightrightarrows X$ with $s(x,g)=x$ and $r(x,g)=g.x$. We denote this by $X \rtimes \Gamma$. A basis $\lbrace U_{i} \rbrace$ for the topology of $X$ lifts to a basis for the topology of $X \rtimes \Gamma$, given by sets $[U_{i},g]:=\lbrace (u,g) \mid u \in U_{i} \rbrace$. 
\end{Ex}

\begin{Ex}
We move now to examples of groupoids coming from uniformly discrete metric spaces of bounded geometry. %Let $X$ be a metric space that is:
%\begin{itemize}
%\item \textit{uniformly discrete}: there exists $c>0$ such that for every pair of distinct points $x,y\in X$, $d(x,y)>c$;
%\item \textit{bounded geometry}: For every $R>0$, there exists a constant $N_{R}$ such that $\vert B_{R}(x) \vert \leqslant N_{R}$ for every $x \in X$. 
%\end{itemize}
We define a groupoid which captures the coarse information associated to $X$. Consider the collection $\mathcal{S}$ of the $R$-neighbourhoods of the diagonal in $X\times X$; that is, for every $R>0$ the set
\begin{equation*}
E_{R}=\lbrace (x,y) \in X \times X \mid d(x,y)\leqslant R \rbrace
\end{equation*}
Let $\mathcal{E}$ be the coarse structure generated by $\mathcal{S}$ as in \cite{MR2007488}; it is called the \textit{metric coarse structure} on $X$. If $X$ is a uniformly discrete metric space of bounded geometry, then this coarse structure is uniformly locally finite, proper and weakly connected -- thus of the type studied by Skandalis, Tu and Yu in \cite{MR1905840}. 

We now define the coarse groupoid following the approach of \cite[Appendix C]{Spakula-Willett}. Let $\beta A$ denote the Stone-\v{C}ech compactification of a set $A$. Set $G(X):=\bigcup_{R>0}\overline{E_{R}}$, where the closure $\overline{E}_{R}$ takes place in $\beta X \times \beta X$ and $G(X)$ has the weak topology coming from the union -- with this topology $G(X)$ is a locally compact, Hausdorff topological space, which becomes a groupoid with the pair groupoid operations from $\beta X\times \beta X$. Another possible approach (for instance that adopted originally in \cite{MR1905840} or in \cite{MR2007488}) is to consider graphs of partial translations on $X$ and form a \textit{groupoid of germs} from this data \cite{MR2419901}. Each approach has value, depending on the particular situation.
\end{Ex}
One advantage of working with groupoids is that they come with many possible \textit{reductions}.

\begin{Def}
A subset of $C\subseteq G^{(0)}$ is said to be \textit{saturated} if for every element $\gamma \in G$ with $s(\gamma) \in C$ we have $r(\gamma) \in C$. For such a subset we can form a subgroupoid of $G$, denoted by $G_{C}$ which has unit space $C$ and $G_{C}^{(2)}=\lbrace (\gamma, \gamma^{'}) \in G^{(2)} \mid s(\gamma),\, r(\gamma)=s(\gamma^{'}), r(\gamma^{'}) \in C \rbrace$. The groupoid $G_{C}$ is called the \textit{reduction of $G$ to $C$}.
\end{Def}

\begin{Rem}\label{rem:imp}
For a uniformly discrete metric space $X$ of bounded geometry there are natural reductions of $G(X)$ that are interesting to consider. It is easy to see that the set $X$ is an open saturated subset of $\beta X$ and in particular this means that the Stone-\v{C}ech boundary $\partial\beta X$ is saturated. We remark additionally that the groupoid $G(X)|_{X}$ is the pair groupoid $X\times X$ (as the coarse structure is weakly connected).
\end{Rem} 

\begin{Def}\label{def:bdry}
The \textit{boundary groupoid} $\partial G(X)$ associated to $X$ is the groupoid reduction $G(X)|_{\partial\beta X}$.
\end{Def}

\subsection{Box spaces as an example}\label{sect:coarse}
Let $\mc X = \lbrace X_{i} \rbrace_{i}$ be a family of finite connected graphs of uniformly bounded vertex degree.

\begin{Def}\label{Def:spaceofgraphs}
The \textit{space of graphs} associated to $\mc X$ is the set $X:=\bigsqcup_{i} X_{i}$, equipped with any metric $d$ that satisfies:
\begin{enumerate}
\item $d|_{X_{i}}$ is the metric coming from the edges of the graph $X_{i}$;
\item $d(X_{i},X_{j})\rightarrow \infty$ as $i+j \rightarrow \infty$.
\end{enumerate}
We remark that any two metrics that satisfy i) and ii) are coarsely equivalent, and thus we need not be more specific about the rates of divergence.
\end{Def}

Natural examples of graph families, and thus spaces of graphs, come from finitely generated residually finite discrete groups. Let $\Gamma=\langle S \rangle$ be finitely generated and residually finite. Then, for any chain (i.e. a nested family of finite index subgroups with trivial intersection) $\mc H = \lbrace H_{i} \rbrace_{i}$ we can consider the Schreier coset graphs:
\begin{equation*}
X_{i}:= \Cay(\Gamma/H_{i}, S).
\end{equation*}

\begin{Rem}
We note that there are various conditions in the literature that one could reasonably put into such a chain of finite index subgroups, for instance asking for each to be normal subgroups, or more generally to separate points from the entire conjugacy class of the subgroup $H_{i}$ (which is called \textit{semi-conjugacy separating} in \cite{finn-sell-wu} and appears first in \cite{Szabo-Wu-Zacharias}), or to ask that the family is \textit{Farber} (that is, for any $g\in \Gamma$, $n_{i}(g) = o(n_{i})$, where $n_{i}$ is the number of conjugates of $H_{i}$ in $\Gamma$ and $n_{i}(g)$ is the number of conjugates of $H_{i}$ containing $g$ \cite{MR1625742,MR2966663}).
\end{Rem}

For simplicity, suppose the chain consists of normal subgroups. Then the space of graphs associated to $\mc X = \lbrace X_{i} \rbrace_{i}$ is called the \textit{box space} of $\Gamma$ with respect to $\mc H$, and denoted by $\square_{\mc H}\Gamma$. 

This construction and the many results concerning it in the literature drive the coarse geometric aspect of this paper. We will focus on the coarse groupoid (and its boundary), to get a better feeling for it in a simpler case than will appear later on.

\begin{Def}\label{Def:MCSI}
Let $\mathcal{S}$ be a family of subsets in $X\times X$.  The family $\mathcal{S}$ \emph{generates $\mathcal{E}$ at infinity} if for every $R>0$ there are finitely many sets $S_1,\dots,S_n\in \mc S$ and a finite subset $F\subset X \times X$ such that
\begin{equation*}
E_{R} \subseteq \left(\bigcup_{k=1}^{n}S_{k}\right)\cup F.
\end{equation*}
\end{Def}

\begin{Rem}
The above definition is equivalent to asking that $\overline{E_{R}}\setminus E_{R} \subseteq \bigcup_{k=1}^{n}\overline{S_{k}}\setminus S_{k}$, where the closure is taken in $\beta X \times \beta X$.
\end{Rem}

If $\Gamma$ is a discrete group acting on $X$, let $E_{g}:= \lbrace (x,x.g)\mid x \in X \rbrace$ be the $g$-diagonal in $X$. We say that the action of $\Gamma$ generates the metric at infinity if the set $\lbrace E_{g} \mid g \in \Gamma \rbrace$ satisfies Definition \ref{Def:MCSI}. 

\begin{Prop}[{\cite[Proposition 2.5]{MR3197659}}]\label{Prop:Crit}
Let $X$ be a uniformly discrete bounded geometry metric space and let $\Gamma$ be a finitely generated discrete group. If $\Gamma$ acts on $X$ so that the induced action on $\beta X$ is free on $\partial \beta X$ and the action generates the metric coarse structure at infinity, then $\partial G(X) \cong \partial \beta X \rtimes \Gamma$. \qed
\end{Prop}

The following example is the basic model we will build on in Section \ref{sect:sofic-coarse-bd-groupoid} for sofic groups.
\begin{Ex}\label{Ex:box}
Let $X=\square_{\mc H}\Gamma$ be the box space of a residually finite group $\Gamma$ with normal chain $\mc H$. 
Then, considering the metric $d$ from Definition \ref{Def:spaceofgraphs} we see that the sets $E_{R}$ decompose as
\begin{equation*}
E_{R} = \bigsqcup_{i} E_{R,i} \sqcup F_{R},
\end{equation*}
where $E_{R,i}$ is the $R$-neighbourhood of the diagonal in $X_{i}$ and $F_{R} = \{(x,y)\mid x \in X_{i}, y \in X_{j}, i\not = j,\; d(x,y)\leqslant R\}$. This observation allows us to reduce to considering the set $E_{R,\infty} = \bigsqcup_{i} E_{R,i} \subset E_{R}$, as these sets have the same Stone-\v{C}ech boundary.

As the group $\Gamma$ is residually finite, each of the $E_{R,i}$ decomposes as $\bigsqcup_{\vert g \vert \leqslant R}E_{g,i}$ when $i$ is sufficiently large -- in particular, $\partial \beta E_{R,\infty} = \bigsqcup_{\vert g \vert \leqslant R} \partial  \beta E_{g}$, and so the group, acting by translations, generates the metric coarse structure at infinity. This action is free at the boundary by residual finiteness of $\Gamma$: for each $g \in \Gamma$ the orbit graph for the action of $g$ on $\square \Gamma$ has degree at most $2$, and thus is at most $3$-coloured by Brookes' theorem. The Stone-\v{C}ech boundaries of each colour set are then permuted by the element $g$ and have empty intersection. Thus Proposition \ref{Prop:Crit} implies that $\partial G(X) \cong \partial \beta X \rtimes \Gamma$. 
\end{Ex}
\subsection{A formal definition of soficity}\label{sect:sofic}

Let us give a formal definition of a sofic group:

\begin{Def}\label{def:sofic}[see \cite[Theorem 3.5]{pestov_hyperlinear_2008}
A group $\Gamma$ is \textit{sofic} if for every finite subset $F \subset \Gamma$ and every $\eps > 0$ there exists a finite set $X$, a map $\sigma \colon \Gamma \rightarrow \Sym(X)$ and a subset $Y\subset X$ with $|Y| \geqslant  (1-\eps)|X|$ such that
\begin{equation*}
\sigma(g)\sigma(h)(y)=\sigma(gh)(y),\quad g,h\in F,\, y\in Y
\end{equation*}
and
\[
\sigma(g)(y) \not = y, \quad g\in F\setminus \lbrace e \rbrace,\, y\in Y.
\]
The map $\sigma$ is said be an $(F,\eps)$\textit{-injective almost action} on the set $X$ if the condition above holds. %We call $Y$ the \textit{sofic core} of $X$.

We note that if $\Gamma$ is sofic, then by fixing a nested sequence of sets $F_{i}$ that exhaust the group, choosing a sequence $\eps_{i} \rightarrow 0$, and letting $X_{i}$ be a set with an $(F_{i},\eps_{i})$-injective almost actions of $\Gamma$, we obtain a sequence of sets together with almost actions of $\Gamma$; such a sequence called a \textit{sofic approximation of} $\Gamma$.
\end{Def}
% The corresponding sequence of the subsets $\lbrace Y_{i} \rbrace$ will be called the sofic core of $\mc X$.

We remark that soficity generalises both being residually finite and being amenable for a group $\Gamma$. We refer the reader to the book \cite{MR3074498} for more details of the permanence properties of sofic groups, and we also note that there is, at time of writing, no group that is known to be non-sofic.

In the remaining part of this section, we will give a more geometric definition of soficity which will allow us to apply coarse geometric methods.

\subsection{Ultralimits and local convergence of graphs}\label{sect:ultralimits}

\begin{Def}
Let $\mc X = \lbrace X_{i} \rbrace_{i}$ be a countable family of finite graphs of bounded degree, $X$ be the space of graphs attached to $\mc X$ and let $\omega \in \partial \beta \mb N$ be a non-principal ultrafilter on $\mathbb{N}$. Let $\underline{x}$ be a sequence of points in $X$, and let $S(\underline{x})$ be the set of all $\underline{y}=(y_{n})_{n}$ such that $\sup_{n}(d(x_{n},y_{n})) < \infty$.
We define a (pseudo-)metric on $S(\underline{x})$ by
\begin{equation*}
d_{\omega}(\underline{y},\underline{z}) = \lim_{\omega}d(y_{n},z_{n})
\end{equation*}
%As the sequences $\underline{y}$ and $\underline{z}$ are at bounded distance from $\underline{x}$ this is a well defined function. 
and the \textit{ultralimit along $\omega$}, denoted $X(\omega,\underline{x})$, to be the canonical quotient metric space obtained from $(S(\underline{x}), d_{\omega})$ by identifying all pairs of points at distance $0$. 
\end{Def}

This notion of ultralimit has a natural description in terms of the coarse boundary groupoid $\mc G:= \partial G(X)$ from the previous section. Let $\eta = \lim_{\omega} x_{n}$ be the point in the Stone-\v{C}ech boundary that corresponds to $\underline{x}$ and $\omega \in \partial\beta \mathbb{N}$. 
 
\begin{Prop}\label{Prop:withRufus}
Let $\mc G_{\eta}$ be the source fibre of $\mc G$ at $\eta \in \partial\beta X$. Equip $\mc G_{\eta}$ with the metric
\begin{equation*}
d_{\eta}((\eta_{1},\eta),(\eta_{2},\eta)) = \inf\lbrace R>0 \mid (\eta_{1},\eta_{2})\in \overline{E_{R}} \rbrace.
\end{equation*}
Then the map $f: X(\omega,\underline{x}) \rightarrow G_{\eta}$ given by $[(y_{n})] \mapsto (\lim\limits_{\omega}y_{n}, \eta)$ is a basepoint preserving isometry.
\end{Prop}
\begin{proof}
For any points $[(y_{n})],[(z_{n})] \in X(\omega,\underline{x})$, we have
\begin{eqnarray*}
d_{\omega}([(y_{n})],[(z_{n})]) & = & \inf\lbrace R>0 \mid \omega(\lbrace n \in \mathbb{N} \mid d(y_{n},z_{n})\leqslant R\rbrace)=1 \rbrace\\
& = & \inf\lbrace R>0 \mid \omega(\lbrace n \in \mathbb{N} \mid (y_{n},z_{n})\in E_{R})=1\rbrace \\
& = & \inf\lbrace R>0 \mid (\lim_{\omega}(y_{n},z_{n})\in \overline{E_{R}}\rbrace \\
& = & \inf\lbrace R>0 \mid (\lim_{\omega}y_{n},\lim_{\omega}z_{n})\in \overline{E_{R}}\rbrace \\
& = & d_{\eta}(\lim_{\omega}y_{n},\lim_{\omega}z_{n}).
\end{eqnarray*}
Hence $f$ is isometric and maps into $\mc G_{\eta}$. It remains to prove that $f$ is surjective.

Let $(\eta',\eta)\in \mc G_{\eta}$. Using the view on $G(X)$ in terms of germs of partial translations as in \cite[Proposition 3.2]{MR1905840} or \cite[Chapter 10]{MR2007488}, we obtain a partial translation $t: A \rightarrow B$ between subsets $A,B \subset X$ such that $\eta \in \overline{A}\subset \beta X$, $\eta' \in \overline{B} \subset \beta X$ and with $\overline{t}(\eta)=\eta'$.% (that is the pair $(\eta',\eta)$ is the unique point with source $\eta$ in the closure of the graph of $t$).
 As $\eta = \lim_{\omega}(x_{n})$, we have that the set $E=\lbrace n \in \mathbb{N} \mid x_{n} \in A \rbrace$ has $\omega$-measure $1$, and therefore we can define another sequence with terms:
\begin{equation*}
y_{n}:=\begin{cases} x_{n} \mbox{ if } n\not\in E\\ t(x_{n}) \mbox{ if } n \in E. \end{cases}
\end{equation*}
As $\eta'$ is the unique point in the closure of the graph of $t$ satisfying $(\eta',\eta) \in \overline{\graph(t)}$, we have that 
\begin{equation*}
(\eta',\eta)=\lim_{\omega}(t(x_{n}),x_{n})=\lim_{\omega}(y_{n},x_{n}),
\end{equation*}
and thus $\eta'=\lim_{\omega}y_{n}$.
\end{proof}

We remark that for a fixed ultrafilter $\eta \in \partial\beta X$ one can always find a sequence $\underline{x}$ tending to infinity and an ultrafilter $\omega \in \partial\beta \mathbb{N}$ such that $\eta = \lim_{\omega} \underline{x}$. There will in general be many such choices, but the above proposition ensures that they will give isometric fibres.

Ideally, we would like to remove the dependence on the base point from this process. The suggested method (say of \cite{MR1873300} or \cite{MR2354165}) is to make this choice \textit{uniformly at random}, and to do this we need a measure on $\partial \beta X$.

Given the sequence of counting measures $\mu_{i}$ on each $X_{i} \in \mc X$ and fixing an ultrafilter $\omega \in \partial\beta \mb N$, we can obtain a measure $\mu$ on the Stone-\v{C}ech boundary of $X$ corresponding to the state
\begin{equation}\label{eq:mu}
\mu(f)=\lim_{\omega} \frac{1}{\vert X_{i}\vert}\sum_{x \in X_{i}} f(x),\quad f\in C(\beta X).
\end{equation}
%and then taking $\mu$ to be the measure obtained by applying the Riesz representation theorem to $\widetilde{\mu}$\footnote{This of course gives many measures as weak cluster points of the counting measures $\mu_{i}$. All of them are equally useful here.}.
Note that $\mu(X)=0$, whence $\mu(\partial\beta X)=1$. Armed with this measure on $\partial\beta X$, we can now formulate a notion of graph convergence:

\begin{Def}
A sequence of graphs  $\mc X$ of bounded degree is said to \textit{Benjamini--Schramm converge} to a graph $Y$ if the set 
\[
\{x = \lim\limits_\omega x_n \in \partial \beta X \mid X(\omega,\underline{x}) \cong (Y,y) \mbox{ for some } y\in Y\}
\]
of ultralimits that are isomorphic as pointed graphs to $Y$ has $\mu$-measure 1.
\end{Def}

A first remark concerning this definition is that the basepoint in $Y$ does not matter if $Y$ is vertex transitive. The second remark we make is that this definition can also be made using \textit{labelled} graphs.

Let $S$ be a finite set of labels. Suppose also that each $X_{i}$ admits an $S$-edge labelling. Then any ultralimit of the sequence $X(\omega, \underline{x})$ also admits an $S$-labelling. In this case, we can ask that $Y$ admits a labelling and that the base point preserving isometries occurring in the definition can be taken as isometries of labelled graphs. 

\begin{Rem}
The traditional formulation of Benjamini--Schramm convergence (found for instance in \cite{MR1873300}) uses converging probabilities of isometry types of balls. It is equivalent to this more topological formulation by realising an ultralimit $X(\omega,\underline{x})$ as a union of balls around $\underline{x}$ and studying how these can be obtained from the sequence $\mc X$ using $\omega$. This works equally well in labelled and non-labelled settings.
\end{Rem}

\begin{Rem}
This Benjamini--Schramm convergence should be thought of as an ``almost everywhere" (in terms of the normalised counting measure) version of the convergence in the space of marked graphs -- if a sequence of bounded degree finite graphs converges there to a fixed graph, then it Benjamini--Schramm converges to that graph -- in fact, the set of measure $1$ will be the entire boundary in that case.
\end{Rem}

The following definition is central to the paper:

\begin{Def}
Let $\Gamma$ be a finitely generated group with a finite generating set $S$. $\Gamma$ is \textit{sofic} if there exists a sequence $\mc X$ of bounded degree, finite $S$-labelled graphs such that $\mc X$ Benjamini--Schramm converges to $(\Cay(G,S),e_{G})$.
\end{Def}

It is equivalent to Definition \ref{def:sofic} by an argument present in \cite[Theorem 5.1]{pestov_hyperlinear_2008}, which constructs the ($S$-labelled) graph structure on the sets $X_{i}$ appearing in Definition \ref{def:sofic} by connecting each $x\in X_i$ with $\sigma_i(s)$ by an edge labelled with $s\in S$; we will always equip $X_i$ coming from a sofic approximation with this graph structure and (slightly abusing notation) also call the resulting sequence $\mc X$ a sofic approximation of $\Gamma$. The following lemma asserts that we can assume these graphs to be connected, which we will always do.

\begin{Lemma}
Let $\Gamma = \langle S \rangle$ be a finitely generated sofic group and let $\mathcal{X}'=\lbrace X'_{i}, \sigma'_{i} \rbrace_{i}$ be a sofic approximation; equip $X'_i$ with the graph structure described above. For each $i$ there is a connected component $X_i\subset X'_i$ and maps $\sigma_i\colon \Gamma\to \Sym(X_i)$ coinciding with $\sigma'_i$ on the generating set $S$ such that $\mathcal{X}=\lbrace X_{i}, \sigma_{i} \rbrace_{i}$ is a sofic approximation with $X_i$. In particular, the graph structure coming from $\mc X$ makes $X_i$ connected.
\end{Lemma}
\begin{proof}
Let $X'_{i,j}$, $j=1,\dots,n_i$ be the connected components of $X'_i$ and let $Y'_i\subseteq X'_i$ be the subsets from Definition \ref{def:sofic}. Increasing $i$ if needed, we may assume without loss of generality that $S\subset F_i$. Observe that
\[
|Y'_i| = \sum_{j=1}^{n_i} |Y'_i\cap X'_{i,j}|\geqslant (1-\eps) |X'_i| = (1-\eps)\sum_{j=1}^{n_i} |X'_{i,j}|.
\]
This implies that there is at least one connected component $X'_{i,j}$ such that $|Y'_i\cap X'_{i,j}|\geqslant (1-\eps) |X'_{i,j}|$; we denote it by $X_i$ and set $Y_i^{(0)} \coloneqq Y'_i\cap X_i$.

Observe that by definition of the graph structure and by preceding construction:
\begin{itemize}
\item the connected components $X'_{i,j}$ are invariant under $\sigma'_i(S)$;
\item $|Y^{(0)}_i| \geqslant (1-\eps)|X_i|$.
\end{itemize}
For $g\in F_i$, we set $Y_{i,g}\coloneqq \{x\in X_i\,|\, \sigma'_i(g)(x)\in X_i\}$. We define $\sigma_i(g) \in \Sym(X_i)$ for $g\in F_i$ by (arbitrarily) extending the partial bijection $\sigma'_i(g)\colon Y_{i,g}\to X_i$ to a permutation $\sigma_i(g) \in \Sym(X_i)$ and we set $\sigma_i(g) = \id_{X_i}$ for $g\not\in F_i$. The above properties guarantee that $\mathcal{X}=\lbrace X_{i}, \sigma_{i} \rbrace_{i}$ is the desired sofic approximation:
\begin{itemize}
\item as $\sigma_i(g)$ coincides with $\sigma'_i(g)$ on the points which remain in $X_i$ under the latter permutation, the set
\[
Y_i \coloneqq \lbrace x \in X_i \mid \forall g,h \in F_i\; \sigma_i(g)\sigma_i(h)(x)=\sigma_i(gh)(x) \mbox{ and } \forall g\in F_i\setminus \lbrace e \rbrace\;\sigma_i(g)(x) \not = x   \rbrace
\]
contains $Y_i^{(0)}$ and therefore satisfies $|Y_i| \geqslant (1-\eps)|X_i|$;
\item $\sigma_i(s) = \sigma'_i(s)$ for all $s\in S$, and therefore the graph structure associated with $\sigma_i$ is the same as the one coming from $\sigma'_i$.  
\end{itemize}
This finishes the proof. 
\end{proof}

\section{The sofic coarse boundary groupoid}\label{sect:sofic-coarse-bd-groupoid}
Let $\Gamma = \langle S \rangle$ be a finitely generated sofic group and $\mc X$ be a sofic approximation of $\Gamma$. The main idea of this paper is that the space of graphs $X$ associated with $\mc X$ can be thought of as a box space for sofic group. In this section we will analyse the boundary groupoid attached with $X$, defined in the previous section. We will also explain how this analysis connects with the sofic core of the sofic approximation.  We remark that being finitely generated by $S$ gives rise to a natural quotient map $\pi_{\Gamma}:F_{S} \rightarrow \Gamma$, where $F_{S}$ is the free group on the letters $S$.

\begin{Def}
Let $\mc G$ be the coarse boundary groupoid associated with the space of graphs $X$ of a sofic approximation $\mathcal{X}=\lbrace X_{i}, \sigma_{i} \rbrace_{i}$ as defined in the previous section. $\mc G$ is called the \emph{sofic coarse boundary groupoid} associated with the sofic approximation $\mc X$. Its base space $\partial\beta X$ is called the \emph{sofic boundary} of $\mc X$.
\end{Def}

\begin{Rem}
For a sofic group $\Gamma$ with a sofic approximation $\mc X$ and the attached space of graphs $X$, for $\mu_{\mc X}$-almost all $\omega \in \beta X$, the range fibre $r^{-1}(\omega)$ is isometric to $\Cay(\Gamma,S)$, as $\mc X$ is a sofic approximation. Let $\delta_{\omega}$ be the Dirac mass at $\omega$ and let $\Ind(\delta_{\omega})$ be the induced representation of $G(X)$ associated with the measure $\delta_{\omega}$ as in \cite{MR2966043}. Then $C^{*}(G(X),\delta_{\omega})$, obtained through the the representation $\Ind(\delta_{\omega})$ of $G(X)$ on $L^{2}(r^{-1}(\omega), \lambda^{\omega})$, is a subalgebra of $C^{*}_{u}(\Gamma)$ \cite[Appendix C]{Spakula-Willett}.
\end{Rem}

As $\mc G$ is a locally compact \'etale groupoid, it can be considered as a Borel groupoid using the natural Borel $\sigma$-algebra obtained from the open subsets of $\mc G$. Our goal in this section is to relate $\mc G$ to an action $\Gamma$, both measurably and topologically. To do this, we introduce an action of $F_{S}$ on $\partial \beta X$. Note that each $X_{i}$ is an $S$-labelled finite graph, with labelled edges constructed using the permutations $\sigma_{i}(s)$. This defines an action of $F_{S}$ on $X_{i}$. We then extend this action continuously to the Stone-\v{C}ech boundary, obtaining an $F_S$-action denoted $\tau$. We remark that when the graphs are regular, it is precisely the action defined in \cite[Lemma 3.26]{MR3197659}.
The action $\tau$ is in general not free, %(unless, for instance the group $\Gamma$ is the free group itself and the quotient map is the identity)
but is still connected with the groupoid $\mc G$.

\begin{Def}
A $\tau$-diagonal on the boundary is a set of the form:
\begin{equation*}
A_{P}:= \lbrace (\omega, \tau(P)(\omega)) \mid \omega \in \partial \beta X \rbrace.
\end{equation*}
for each $P \in F_{S}$.
\end{Def}

\begin{Prop}\label{Prop:coarseER}
$\mc G$ is isomorphic to the orbit equivalence relation $\mc R_{\tau}$ of the action $\tau: F_{S} \rightarrow \Homeo(\partial\beta X)$, where this equivalence relation is given the  weak topology generated by the clopen sets $\lbrace A_{P} \rbrace_{P \in F_{S}}$.
\end{Prop}
\begin{proof}
We check that, for each $n \in \mathbb{N}$, the sets $\partial E_{n}$ and $\bigcup_{\vert P \vert \leq n} A_{P}$ are equal. We first observe that if $\gamma \in \partial E_{n}$ then there is a net of pairs $((x_{\lambda},y_{\lambda}))_{\lambda}$ with limit $\gamma$, and $d(x_{\lambda},y_{\lambda})\leqslant n$ on a convergent subnet. 

However, as the distance here is natural edge metric on a graph, to be at distance of at most $n$ means that $x_{\lambda}$ and $y_{\lambda}$ are connected by an $S$-labelled path of length of most $n$. From this we conclude that the $F_{S}$-action by the concatenation of the labels will map $x_{\lambda}$ to $y_{\lambda}$. 

To see the reverse inclusion, we observe that anything belonging to at least one of the $A_{P}$'s must be a limit of a net of pairs of the form $(x_{\lambda}, \tau(P)(x_{\lambda}))$. Therefore this net consists of pairs whose distances are bounded precisely by the length of $P$, which was supposed less than $n$.
\end{proof}

We now return to $\Gamma$. For each $g \in \Gamma$, the map $\sigma(g)$ defined by performing $\sigma_{i}(g)$ in each graph $X_{i}$ defines a bijection of $X$ to itself. Extending these maps continuously gives us a collection of homeomorphisms $\sigma(g)$ on $\beta X$. We remark that this gives a map $\Gamma \rightarrow \Homeo(\partial\beta X)$, which is in general \textit{not} a homomorphism of groups, but it is quite close to a homomorphism when we make use of the fact that the soficity of $\Gamma$ is being witnessed by $\mc X$.

Let $Y\subset X$ be the the disjoint union of each $Y_{i}$ coming from Definition \ref{def:sofic}. As the sets $Y_{i}^{c}$ are at most $\mu_{i}$-measure $\eps_{i}$ (and tending to $0$) we have that $\mu(\overline{Y}) = 1$, where $\mu$ is a probability measure on $\partial\beta X$ defined in \eqref{eq:mu}. For any element $\omega \in \partial Y$, the maps $\sigma(g)\sigma(h)$ and $\sigma(gh)$ coincide, and thus the map $\sigma$ is a homomorphism of groups after throwing out a set of measure $0$ in $\partial\beta X$. In particular, this is an example of a ``near action'' of $\Gamma$ in the sense of \cite{MR2191233}.

This is not yet useful topologically, but we can still make the following definition:

\begin{Def}
The $\sigma$-diagonals in $\partial\beta X \times \partial \beta X$ are sets of the form:
\begin{equation*}
E_{g}:= \lbrace (x,\sigma(g)x) \mid x \in X \rbrace,
\end{equation*}
for $g \in \Gamma$.
\end{Def}

Now we relate the equivalence relation $\mc R_{\tau}$ to the  $\Gamma$-near action on $\partial \beta X$ by finding an $F_{S}$-invariant subset of $\partial \beta X$ on which the free group action really agrees with the $\Gamma$-near action.

\begin{Def}
The set
\begin{equation*}
Z:= \bigcap_{g \in \Gamma}\sigma(g)(\partial Y)
\end{equation*}
is called the \emph{core} of the sofic boundary $\partial\beta X$. It depends on the choice of the subsets $Y_i\subset X_i$ satisfying the conditions of Definition \ref{def:sofic}. %Taking $Y_i$ to be maximal subsets of $X_i$ with respect to this condition, we obtain a set $Z_\max\subset \partial\beta Y$ called the \emph{maximal core}. 
\end{Def}
As $\partial Y$ is clopen and the maps $\sigma(g)$ are all homeomorphisms, the core $Z$ is a closed subset of $\partial \beta X$ that is invariant under the maps $\sigma(g)$. Using de Morgan's law, it's clear that $\mu(Z)=1$; in particular the core is not empty.

% \begin{Rem}
% There is a \Cs-algebraic construction yielding the sofic core. As $C(\partial\beta X) = \ell^\infty(X)$, the maps $\sigma(g)$ define automorphisms of $\ell^\infty(X)$ which we denote $\sigma_g$. Now, consider the ideal
% \[
% I \coloneqq \langle \sigma_{gh}(f) -\sigma_g\circ \sigma_h(f),\, f\in \ell^\infty(X)\rangle
% \]
% and let $Z$ be the unique closed subspace of $X$ such that $I = \ker (r_Z\colon \ell^\infty(X) = C(\partial\beta X) \to C(Z))$, where $r_Z$ is the natural restriction map.
% \end{Rem}

 For $K\subset \partial \beta X \times \partial \beta X$, we denote by $K^{Z}$ the restriction $K \cap (Z \times Z)$.

\begin{Lemma}\label{Lem:freevssofic}
We have the following compatibility between the action of $F_{S}$ and the action of $\Gamma$ on $Z$:
\begin{enumerate}
\item For $g\not = h \in \Gamma$, we have that $\partial E_{g}^{Z}\cap \partial E_{h}^{Z} = \varnothing$.
\item $\Stab_{F_{S}}(Z) = \ker(\pi_{\Gamma}:F_{S}\rightarrow \Gamma)$;
\item If $\pi_{\Gamma}(P)=\pi_{\Gamma}(Q)$ then $A_{P}^{Z}=A_{Q}^{Z}$.
\end{enumerate}
\end{Lemma}
\begin{proof}
For i), let $(\omega, \sigma(g)(\omega)) = (\omega,\sigma(h)(\omega)) \in \partial E_{g}^{Z}\cap \partial E_{h}^{Z}$. Thus, $\omega = \sigma(g)^{-1}\sigma(h)(\omega)$. As $Z \subset \partial Y$, we have that $\omega = \sigma(g^{-1}h)(\omega)$, however this can only happen if $g^{-1}h = e$.

The proofs of the remaining points follow directly from a key observation that comes from the definition of $Z$: if $w = a_{s_{1}}\cdots a_{s_{n}} \in F_{S}$, then $\tau(w)(\omega)=\sigma(s_{1})\cdots\sigma(s_{n})(\omega)=\sigma(\pi_{\Gamma}(w))(\omega)$ for every $\omega \in Z$. ii) and iii) are now deduced by elementary calculations using this observation.
\end{proof}

We conclude that the set $Z$ is a closed subset which is invariant under the equivalence relation $\mc R_{\tau}$, and thus under $\mc G$. In fact, combining with the arguments in the proof of Proposition \ref{Prop:coarseER}, we can observe:

\begin{Lemma}\label{Lem:freevssofic-homeo}
There is a homeomophism $\partial E_{n}^{Z} = \bigsqcup_{\vert g \vert \leqslant n} \partial E_{g}^{Z}$, given explicitly by the map
\[
\Theta: \partial E_{n}^{Z} \rightarrow \bigsqcup_{\vert g \vert \leqslant n}\partial E_{g}^{Z},
\]
\[
\gamma \mapsto (s(\gamma), \pi_{\Gamma}(P)(s(\gamma))).
\]
\qed
\end{Lemma}

The main result of this section is the following:

\begin{Thm}\label{Thm:clevertrick}
The reduction groupoid $\mc G|_{Z}$ and the transformation groupoid $Z \rtimes \Gamma$ are topologically isomorphic.
\end{Thm}
\begin{proof}
The technique of the proof is similar to that of Proposition \ref{Prop:coarseER}.
As $\mc G|_{Z} = \bigcup_{n} \partial E_{n}^{Z}$, and $Z \rtimes \Gamma$ is the disjoint union $\bigsqcup_{g \in \Gamma} \partial E_{g}^{Z}$, we obtain a map $\Theta: \mc G|_{Z} \rightarrow Z\rtimes \Gamma$ using the (obviously compatible) map from Lemma \ref{Lem:freevssofic-homeo}. It remains to see that it is both a homeomorphism and a homomorphism of groupoids.

We observe that: 
\begin{enumerate}
\item both groupoids have a basis of topology given by clopen slices \cite[Proposition 4.1]{MR2644910};
\item as $\mc G$ has the weak topology, it is sufficient to consider slices contained in $\overline{E_{n}}$, i.e we can consider slices $U \subset \partial E_{n}^{Z}$ when working with $\mc G|_{Z}$;
\item slices for $Z\rtimes \Gamma$ are of the form $(U,g):=\lbrace (\omega, \sigma(g)\omega) \mid \omega \in U\}$ for some clopen $U \subset Z$.
\end{enumerate}
Given a slice $U \subset \mc G|_{Z}$ contained in some $\partial E_{n}^{Z}$, we can see that $\Theta(U)$, by Lemma \ref{Lem:freevssofic} iv), is contained within a finite disjoint union of clopen sets $\partial E_{g}^{Z}$. This means, in particular, that $\Theta(U)=\bigsqcup_{g}(U,g)$, which are open and disjoint. A similar argument proves that the map $\Theta$ is continuous. 

To complete the proof we must show that the map is a homomorphism. This, however, follows from Lemma \ref{Lem:freevssofic} ii) and the fact the map $\pi_{\Gamma}:F_{S} \rightarrow \Gamma$ is a group homomorphism.
\end{proof}

Recall that the measure $\mu$ is naturally extended to a Borel measure $\nu:=\mu\circ \lambda$ on $\mc G|_{Z}$, defined by:
\begin{equation*}
\int_{\gamma \in \mc G} f d\nu = \int_{x \in \partial\beta X} \left(\sum_{s(\gamma)=x} f(\gamma)\right) d\mu(x)
\end{equation*}
for every Borel measurable function $f$ on $\mc G|_{Z}$.

\begin{Cor}
The measure $\nu = \mu \circ \lambda$ is invariant for $\mc G|_{Z}$ (and thus for $\mc G$). 
\end{Cor}
\begin{proof}
We compute:
\begin{equation*}
\int_{\gamma \in \mc G|_{Z}} f d\nu = \sum_{g \in \Gamma}\int_{\gamma \in \partial E_{g}^{Z}} f d\nu.
\end{equation*}
We now analyse the last integral under the map $\gamma \mapsto \gamma^{-1}$, where it transforms to:
\begin{equation*}
\int_{\gamma^{-1} \in \partial E_{g}^{Z}} f d\nu = \int_{x \in Z}\sum_{\substack{s(\gamma^{-1})=x\\ \gamma^{-1}\in \partial E_{g}^{Z}}} f(\gamma^{-1}) d\mu(x).
\end{equation*}
The conditions on the integrand here are equivalent to the statement that $\gamma \in \partial\Delta_{g^{-1}}^{Z}$ and that $s(\gamma)=\sigma(g)(x)$. As $\mu$ and $Z$ are both invariant under $\sigma(g)$, performing a change of variables $x \mapsto \sigma(g)^{-1}(x)$ we see that this last integral is equal to:
\begin{equation*}
\int_{x \in Z}\sum_{\substack{s(\gamma)=x\\ \gamma \in \partial E_{g^{-1}}^{Z}}} f(\gamma) d\mu(x)=\int_{\gamma \in \partial E_{g^{-1}}^{Z}} f d\nu.
\end{equation*}
However, as we are summing over the group $\Gamma$, this completes the proof.
\end{proof}

Thus $(\mc G|_{Z}, \nu)$ is a measured groupoid and the topological isomorphism of Theorem \ref{Thm:clevertrick} gives us an isomorphism of measured groupoids $(\mc G|_{Z}, \nu) \cong (Z,\mu)\rtimes \Gamma$. Thus, if we extend the action of $\Gamma$ on $\partial\beta X$ by letting every element of $\Gamma$ act by the identity on the complement of $Z$, we obtain an almost everywhere isomorphism\footnote{This is just an isomorphism in parts of the measured groupoid literature, cf. \cite{MR3130315}.} as in \cite{MR665021} for $\mc G$ and $\partial\beta X \rtimes \Gamma$:

\begin{Thm}\label{thm:ae-iso}
The measured groupoids $(\mc G,\nu)$ and $(\partial\beta X, \mu)\rtimes \Gamma$ (where each element of $\Gamma$ is defined to act by the identity transformation on the complement of $Z$) are almost everywhere isomorphic as Borel measured groupoids.
\end{Thm}
\begin{proof}
The map defined in the proof of Theorem \ref{Thm:clevertrick} is a well defined groupoid homomorphism of topological groupoids, but the set of elements in $\mc G$ for which this map is not well defined have measure $0$; this is precisely the definition of an almost everywhere isomorphism: just map the elements $\gamma = (\omega, \omega^{'}) \in \mc G|_{Z^{c}}$ to any pair $(\omega, \tau(P_{\gamma}))$ and notice that the homomorphism rule will hold almost everywhere for the appropriate measure on $\mc G$.
\end{proof}

\begin{Rem}
In the purely measurable setting, given a sofic approximation $\mc X$ and an ultrafilter $\omega\in\partial\beta \mb N$, one can naturally define the ultraproduct measure space
\[
\prod_{i\to\omega} (X_i,\mu_i)
\]
which will carry a natural $\Gamma$-action: viewing the sofic approximation $\sigma$ as an embedding of $\Gamma$ into the ultraproduct of permutation groups,
\[
\sigma\colon \Gamma\hookrightarrow \prod_{i\to \omega} \Sym(X_i),
\] 
one uses natural embeddings $\Sym(X_i)\hookrightarrow \mb M_{|X_i|}(\mb C)$ as permutation matrices to obtain a unitary representation
\[
\sigma\colon \Gamma\hookrightarrow \mc U\left(\prod_{i\to \omega} (M_{|X_i|}(\mb C),\tr_i)\right),
\]
where $\tr_i$ denotes the normalized trace. As permutation matrices normalize the subalgebra of diagonal matrices $A_i \subset M_{|X_i|}(\mb C)$, we obtain a natural action of $\Gamma$ on the ultraproduct
\[
\Gamma \curvearrowright \prod_{i\to \omega} (A_i,\tr_i),
\]
and this latter ultraproduct is by construction isomorphic to
\[
\prod_{i\to \omega} (A_i,\tr_i) \cong \prod_{i\to \omega} (\ell^\infty(X_i),\mu_i) \cong L^\infty\left(\prod_{i\to\omega} (X_i,\mu_i)\right)
\]
On the other hand, by definition of the ultraproduct
\[
\prod_{i\to \omega} (\ell^\infty(X_i),\mu_i) \cong \ell^\infty(X)/\{f\in\ell^\infty(X)\,|\, \lim\limits_{i\to\omega} \mu_i(f^* f) = 0\} \cong L^\infty(\partial\beta X,\mu).
\]
Therefore measure theoretically our construction yields nothing but the ultraproduct measure space naturally associated with the sofic approximation.
\end{Rem}

\begin{Rem}
The results in this section should be thought of as an ``almost everywhere" version of Example \ref{Ex:box}, where the set $Z$ should be considered as the appropriate boundary set to attach to the space of graphs $X$ of a sofic approximation $\mc X$.
\end{Rem}

\section{From sofic approximations to analytic properties of the group}\label{sect:main-theorems}
In this section we prove the results announced in Theorem \ref{thm:intro}, and we recall the necessary definitions (or references) of the coarse geometric and analytic properties that we need to keep this paper approximately self contained.
 
\subsection{Amenability}

Let $X$ be a uniformly discrete metric space of bounded geometry. We begin with a few definitions concerning $X$:
\begin{Def}
$X$ is \textit{amenable} if for every $R>0,\eps>0$ there exists a finite set $F\subset X$ such that 
\begin{equation*}
\frac{\vert \partial_{R}F \vert}{\vert F \vert} < \eps,
\end{equation*}
where $\partial_{R}F$ is the $R$-boundary of $F$, that is the set of points in the $R$-neighbourhood of $F$ that do not themselves belong to $F$.
\end{Def}

Equivalent to this metric definition is a functional one:

\begin{Def}
 $X$ is $(R,\eps)$-amenable if there exists a norm one probability measure $\phi$ on $X$ such that:
\begin{equation*}
\sum_{(x,y)\in E_{R}}\vert \phi(x) - \phi(y) \vert \leqslant \eps.
\end{equation*}
\end{Def}

A space $X$ is \textit{amenable} if it is $(R,\eps)$-amenable for every $R>0$, $\eps>0$ \cite{MR1145337}.

This leads nicely to a functional definition of \textit{property A}, a coarse notion of amenability introduced by Yu in \cite{MR1728880}, which is heavily studied in the literature. For a comprehensive survey on what is known about property A, see \cite{MR2562146}.
\begin{Def}
 $X$ has \textit{Property A} if for every $R>0, \eps>0$, there exists an $S>0$ and a function $\eta \colon X \rightarrow \Prob(X)$, written $x \mapsto \eta_{x}$ with the following properties:
\begin{enumerate}
\item each $\eta_{x}$ is supported in a ball of radius at most $S$ around $x$;
\item for any pair $(x,y)\in E_{R}$, we have: $\Vert \eta_{x} - \eta_{y} \Vert \leqslant \eps$.
\end{enumerate}
Condition ii) for $\eta$ is known as being $(R,\eps)$-variation.
\end{Def}
For \textit{families} of metric spaces, we can study \textit{uniform} properties of the family. In this context, a family $\mc X=\lbrace X_{\alpha} \rbrace_{\alpha}$ has property A uniformly if, for every $R>0,\eps>0$ and there is an $S>0$ independent of $\alpha$ such that $X_{\alpha}$ satisfies conditions in the definition of propety A for parameters $R,\eps, S$. 

\begin{Ex}
For families of metric spaces, we know the following:
\begin{enumerate}
\item Any sequence of finite graphs $\lbrace X_{i} \rbrace_{i}$ with degree bounded below by $3$, above uniformly and girth tending to $\infty$,  does not have property A uniformly, where girth is the length of the shortest simple cycle \cite{MR2831267};
\item Any box space of any residually finite amenable group is property A (in fact, this characterises amenability for a residually finite group) \cite[Chapter 11]{MR2007488}.
\end{enumerate}
\end{Ex}

Here is the amenability part of the main result of this paper.

\begin{Thm}\label{Thm:amenable}
Let $\Gamma$ be a sofic group and let $\mc X$ be the a sofic approximation of $\Gamma$. If $\mc X$ has property A uniformly, then $\Gamma$ is amenable. %Moreover, if $Y$ is the sofic core of the space of graphs $X$ of $\mc X$, then $\Gamma$ is amenable if and only if $Y$ has property $A$.
\end{Thm}
\begin{proof}
Suppose the space of graphs $X$ associated with $\mc X$ is property A. Then the full coarse groupoid -- and thus $\mc G$, which is a closed reduction -- is topologically amenable as a groupoid \cite{MR1905840}. Applying this closed reduction fact again, $Z \rtimes \Gamma$ is therefore topologically amenable -- but since $Z$ has a $\Gamma$-invariant probability measure, this can happen if and only if $\Gamma$ is amenable \cite[Example 2.7.(3)]{a-d-amenability-exactness}. 
\end{proof} 

\subsection{Amenable limits}

As a basic application of the ideas from Section \ref{sect:ultralimits}, we also give an answer to the following natural question: given a graph sequence with property A, can one use the measure $\mu$ to tell ``how many'' ultralimits are amenable as metric spaces?

Let $A_{\mathrm{amen}}$ denote the set of ultralimits of a graph sequence $\mc X$ that are amenable as metric spaces.

\begin{Prop}\label{prop:withRufus}\quad
\begin{enumerate}
\item If $\mc X=\lbrace X_{i} \rbrace_{i}$ is a family of finite graphs with bounded degree that has property A uniformly, then there exists an ultralimit $X(\omega, \underline{x})$ that is $(R,\eps)$-amenable;
\item If $\mc X$ has property A and Benjamini--Schramm converges to a graph $X$, then $\mu(A_{\mathrm{amen}})\in \{0,1\}$;
\item For every $q \in \mathbb{Q}\cap [0,1]$ there is a sequence of finite graphs $\mc X$ of bounded degree that have $\mu(A_{\mathrm{amen}})=q$.
\end{enumerate}
\end{Prop}
\begin{proof}
For i): as $\mc X$ is property A uniformly, for each $R,\eps>0$ we can find an $S>0$ (independent of $i$) and a function, for each $i$:
\begin{equation*}
\eta\colon X_{i} \rightarrow \Prob(X_{i}),
\end{equation*}
satisfying:
\begin{itemize}
\item each $\eta_{x}$ is supported in a ball of radius at most $S$ around $x$;
\item for any pair $(x,y)\in E_{R}$, we have: $\Vert \eta_{x} - \eta_{y} \Vert \leqslant \eps/ N_{R}$,
\end{itemize}
where $N_{R}$ is the uniform upper bound on the cardinality of a ball of radius $R$ in $X_{i}$. 

We now unpack the latter point (and using $\Vert \eta_{x} \Vert =1$) into:
\begin{equation*}
\sum_{z \in X_{i}} \vert \eta_{x}(z) - \eta_{y}(z) \vert \leqslant \frac{\eps}{N_{R}}\sum_{z \in X_{i}}\vert \eta_{x}(z) \vert.
\end{equation*}
Fixing $x \in X_{i}$ and summing over the ball of radius $R$ around $x$ gives:
\begin{equation*}
\sum_{z \in X_{i}} \sum_{y \in B_{R}(x)} \vert \eta_{x}(z) - \eta_{y}(z) \vert \leqslant \eps\sum_{z \in X_{i}}\vert \eta_{x}(z) \vert.
\end{equation*}
Now summing over all possible $x \in X_{i}$, we obtain
\begin{equation*}
\sum_{z \in X_{i}} \sum_{(x,y) \in E_{R}} \vert \eta_{x}(z) - \eta_{y}(z) \vert \leqslant \eps\sum_{z \in X_{i}}\sum_{x \in X_{i}}\vert \eta_{x}(z) \vert.
\end{equation*}
It follows from this that there must be some $z_{i} \in X_{i}$ such that:
\begin{equation*}
\sum_{(x,y) \in E_{R}} \vert \eta_{x}(z) - \eta_{y}(z) \vert \leqslant \eps\sum_{x \in X_{i}}\vert \eta_{x}(z) \vert.
\end{equation*}
This lets us define $\phi: X_{i} \rightarrow [0,1]$ by $\phi(x) = \eta_{x}(z_{i})$, and then by the above we deduce:
\begin{equation*}
\sum_{(x,y) \in E_{R}} \vert \phi(x) - \phi(y) \vert \leqslant \eps \Vert \phi \Vert_{1}.
\end{equation*}
and as $\eta_{x}$ is supported in a ball of radius $S$ for each $x$, $\phi$ also is supported in a ball of radius $S$.

Repeating this for each $X_{i}$ and renormalising, we see that for every $R>0,\eps>0$ there exists $S>0$ such that for every $i \in \mathbb{N}$ there is an $z_{i} \in X_{i}$ and a function $\phi_{i}:X_{i} \rightarrow [0,1]$ supported in the ball of radius $S$ around $z_{i}$ such that:
\begin{equation*}
\sum_{(x,y) \in E_{R}} \vert \phi_{i}(x) - \phi_{i}(y) \vert \leqslant \eps.
\end{equation*}
Now take $\underline{z}=(z_{i})_{i}$ and fix any nonprincipal ultrafilter $\omega \in \partial\beta\mathbb{N}$. We claim that the ultralimit $X(\omega,\underline{z})$ is $(R,\eps)$-amenable. Indeed, if we let $B = B_{R+S}(\underline{z})$ in $X(\omega,\underline{z})$, then the set:
\begin{equation*}
E = \lbrace i \in \mathbb{N} \mid B_{R+S}(x_{i}) \mbox{ is isometric to } B \rbrace
\end{equation*}
has $\omega$-measure $1$.

Now, for each $i \in E$ we can use a fixed isometry to transplant $\phi_{i}$ onto the set $B$. We note that these new transplanted functions also satisfy:
\begin{equation*}
\sum_{(x,y) \in E^{X(\omega,\underline{z})}_{R}} \vert \phi_{i}(x) - \phi_{i}(y) \vert \leqslant \eps.
\end{equation*}
As $B$ is bounded, we can now take the ultralimit $\phi = \lim_{\omega}\phi_{i}$, which now clearly satisfies:
\begin{equation*}
\sum_{(x,y) \in E^{X(\omega,\underline{z})}_{R}} \vert \phi(x) - \phi(y) \vert \leqslant \eps.
\end{equation*}

For ii), observe that a graph family $\mc X$ converges to a graph $X$ locally implies that $\mu$-almost all $X(\omega,\underline{x})$ are isometric to $X$, that is we can find a base point $x \in X$ and a basepoint preserving isometry $X(\omega,\underline{x}) \rightarrow (X,x)$ for almost all admissible sequences $\underline{x}$.

Running the proof of i) sequentially for the sequence $(R_n,\eps_n) = (n,\frac{1}{n})$, we construct a family of ultralimits denoted by $Y_{n}$. Now, either $Y_{n}$ is isometric to $X$ for arbitrarily large $n$, or it isn't -- and the first case gives us that $X$ is amenable (as it's $(R,\eps)$-amenable for all $R,\eps>0$). To complete the proof, notice that because of the local convergence, the second case happens for a set of possible admissible sequences of $\mu$-measure $0$.

For iii): fix $q = \frac{a}{b} \in \mathbb{Q}$. Consider the graph family $\mc X = \lbrace X_{i} \rbrace_{i}$ with
\begin{equation*}
X_{i}=\bigsqcup_{k=1}^{a} Y_{i} \sqcup \bigsqcup_{k=a+1}^{b} Z_{i},
\end{equation*}
where $Y_{i}$ is a cycle of length at $i$ and $Z_{i}$ is a family of bounded degree graphs with all vertices of degree at least three and girth at least $i$. Let $X$ be the space of graphs attached with $\mc X$, and let $Y$ and $Z$ be the spaces of graphs attached with the sequences $\mc Y=\lbrace Y_{i} \rbrace_{i}$, $\mc Z=\lbrace Z_{i} \rbrace_{i}$ respectively. Then the boundary $\partial\beta X$, by definition, splits into $\bigsqcup_{k=1}^{a}\partial \beta Y \sqcup \bigsqcup_{k=a+1}^{b} \partial\beta Z$, and thus $\mu\left(\bigsqcup_{k=1}^{a}\partial \beta Y\right)=q$. So for the first part of the claim, it is enough to see that $A_{\mathrm{amen}}=\bigsqcup_{k=1}^{a}\partial \beta Y$. This is clear, however, as any ultralimit of the sequence $Z_i$ is an infinite tree with all vertices
of degree at least three, which is certainly not amenable (this proves $A_{\mathrm{amen}} \subset \bigsqcup_{k=1}^{a}\partial \beta Y$). For the other inclusion, notice that any ultralimit attached the sequence $\mc Y$ is a copy of the integer bi-infinite ray -- this is certainly amenable as a metric space (using the Følner argument for the integers).
\end{proof}

\subsection{a-T-menability}

The following is a compression of definitions taken from \cite{MR1703305} and \cite{MR3142029}.

\begin{Def} 
Let $G$ be a groupoid.
\begin{itemize}
\item A (real) conditionally negative definite function on $G$ is a function $\psi\colon G \rightarrow \mathbb{R}$ such that:
\begin{enumerate}
\item $\psi(x) = 0$ for every $x \in G^{(0)}$;
\item $\psi(g) = \psi(g^{-1})$ for every $g \in G$;
\item For every $x \in G^{(0)}$ , every $g_{1},...,g_{n}\in G^{x}$, and all real numbers $\lambda_{1},...,\lambda_{n}$ with $\sum_{i=1}^{n}\lambda_{i} = 0$ we have:
\begin{equation*}
\sum_{i,j}\lambda_{i}\lambda_{j}\psi(g_{i}^{-1}g_{j}) \leqslant 0
\end{equation*}
\end{enumerate}
\item A locally compact, Hausdorff groupoid $G$ is \textit{a-T-menable} if there exists a proper, continuous, conditionally negative definite function $\psi\colon G \rightarrow \mathbb{R}$. This definition applies to groups: a group $\Gamma$ is a-T-menable if is satisfies ii).
\item A Borel groupoid $(G,\nu)$ is \textit{a-T-menable} if there exists a proper, Borel, conditionally negative definite function $G \rightarrow \mathbb{R}$. In this context, properness means that $\nu(\lbrace g \in G \mid \psi(g) \leqslant c \rbrace )< \infty$ for every $C>0$.
\end{itemize}

If $G$ is locally compact, Hausdorff, topologically a-T-menable groupoid, then the associated Borel groupoid $(G,\nu_{\mu})$ is a-T-menable in the sense of iii) for any quasi-invariant measure $\mu$ on $G^{(0)}$. It's also transparent that topological a-T-menability passes to closed subgroupoids.
\end{Def}

Related to this are various notions of a coarse embedding for a metric space $X$.

\begin{Def}
A metric space $X$ \textit{coarsely embeds} into Hilbert space $H$ if there exist maps $f:X \rightarrow H$,  and non-decreasing $\rho_{1},\rho_{2}:\mathbb{R}_{+} \rightarrow \mathbb{R}$ such that:
\begin{enumerate}
\item for every $x,y \in X$, $\rho_{1}(d(x,y)) \leqslant \Vert f(x) - f(y) \Vert \leqslant \rho_{2}(d(x,y))$;
\item for each $i$, we have $\lim_{r \rightarrow \infty}\rho_{i}(r) = +\infty$.
\end{enumerate}
\end{Def}

The connection with groupoids here is that a result of \cite{MR1905840}, which states that $X$ coarsely embeds into Hilbert space if and only if $G(X)$ is topologically a-T-menable. In \cite{MR3346926}, Willett introduced a property sufficient for the a-T-menability of the boundary groupoid associated with a sequence of bounded degree graphs:

\begin{Def}
Let $\mc X = \lbrace X_{i} \rbrace_{i}$ be a sequence of finite graphs of bounded degree. Then the sequence $\mc X$ \textit{asymptotically (coarsely) embeds into Hilbert space} if there exist non-decreasing control functions $\rho_{1},\rho_{2}:\mathbb{R}_{+} \rightarrow \mathbb{R}$ and symmetric, normalised kernels:
\begin{equation*}
K_{i}: X_{i}\times X_{i} \rightarrow \mathbb{R},
\end{equation*}
and a sequence of non-negative real numbers $(R_{i})_{i}$ tending to infinity satisfying:
\begin{enumerate}
\item for all $i$, and all $x,y \in X_{i}$:
\begin{equation*}
\rho_{1}(d(x,y)) \leqslant K_{i}(x,y) \leqslant \rho_{2}(d(x,y));
\end{equation*}
\item for any $i$ and any subset $\lbrace x_{1},...,x_{n}\rbrace \subset X_{i}$ of diameter at most $R_{i}$, and any collection of real numbers $\lambda_{1},...,\lambda_{n}$ with $\sum_{i} \lambda_{i} =0 $ we have:
\begin{equation*}
\sum_{i,j}\lambda_{1}\lambda_{2}K_{i}(x_{i},x_{j}) \leqslant 0.
\end{equation*}
\end{enumerate}
\end{Def}

The key point here is the parameter family $(R_{i})_{i}$. If this sequence grows faster than the sequence of diameters, then the family $\mc X$ is coarsely embeddable into Hilbert space (uniformly in $i$). However, this might grow slower than the diameter as is the case when the space $X$ \textit{fibred coarsely embeds into Hilbert space} but does not coarsely embed into Hilbert space. The following is \cite[Lemma 5.3]{MR3346926}, which is proved using the techniques of \cite{MR3266245}:

\begin{Prop}
If $\mc X$ is an asymptotically coarsely embeddable family of finite graphs of bounded degree, then the boundary groupoid $\mc G$ of the associated space of graphs $X$ is topologically a-T-menable.\qed
\end{Prop}

Let $\mc G$ be the coarse boundary groupoid of the graphs obtained from the sofic approximation and $Z \subset \partial\beta X$ be a core of the sofic boundary.

\begin{Prop}\label{Prop:Multipliers}
If $\mc G |_{Z}$ is a-T-menable, then $\Gamma$ is a-T-menable.
\end{Prop}
\begin{proof}
As $\mc G|_Z \cong Z\rtimes \Gamma$ and carries an invariant measure, in view of \cite[Corollary 5.11]{MR3138486} it is enough to prove that the action of $\Gamma$ on $Z$ is a-T-menable in the sense of \cite[Definition 5.5]{MR3138486}; this, however, immediately follows from a-T-menability of $\mc G|_Z \cong Z\rtimes \Gamma$.

\end{proof}

\begin{Thm}\label{Thm:a-t-men}
If $\Gamma$ is a sofic group admitting a sofic approximation $\mathcal{X}$ that asymptotically embeds into Hilbert space. Then $\Gamma$ is a-T-menable.
\end{Thm}
\begin{proof}
As $X$ asymptotically coarsely embeds into Hilbert space, the groupoid $\mc G$ is topologically a-T-menable. As $\mc G|_{Z}$ is closed, it also topologically a-T-menable. The result now follows from Proposition \ref{Prop:Multipliers}.
\end{proof}

%The authors would like to remark that the proof given in Proposition \ref{Prop:Multipliers} uses the \textit{measured groupoid} structure attached to $Z \rtimes \Gamma$. If the reader would prefer, they can instead prove this result using the \textit{topological groupoid} structure by appealing instead to Corollary 5.11 of \cite{MR3138486}.

\subsection{Property (T)}

\begin{Def}
A finitely generated discrete group $\Gamma=\langle S \rangle$ has \textit{property (T)} if for any unitary representation $\pi:\Gamma \rightarrow \mc U(H)$ that has \textit{almost invariant vectors} has an invariant vector. Here, a vector $v \in H$ is $\eps$-invariant If
\begin{equation*}
\max_{s\in S} \Vert \pi(s)v - v \Vert \leqslant \eps,
\end{equation*}
and $\pi$ has \textit{almost invariant vectors} if for every $\eps>0$ there is a $\eps$-invariant vector.
\end{Def}

Given a uniformly discrete metric space $X$ of bounded geometry, there is a way to associate a $C^{\ast}$-algebra to $X$ that bridges operator algebraic properties with coarse geometric properties. Let $\ell^{2}(X)$ be the complex Hilbert space spanned by Dirac functions $\delta_{x}$ for each point $x\in X$. Any bounded linear operator $T\in \mb{B}(\ell^{2}(X))$ can be uniquely represented as a matrix $(T_{x,y})$ indexed by $X \times X$ where the entries are defined by $T_{x,y}=\langle T\delta_{x},\delta_{y} \rangle$. 

For $T \in \mb{B}(\ell^{2}(X))$ we can define the propagation of $T$  by the formula:
\begin{equation*}
\mathop\mathrm{Propagation}(T):= \sup \lbrace d(x,y) \mid T_{x,y} \not = 0\rbrace. 
\end{equation*}

\begin{Def}The $\ast$-subalgebra of $\mc B(\ell^2 X)$ consisting of operators with finite propagation is denoted $\mathbb{C}[X]$. The closure of $\mathbb{C}[X]$ in the operator norm of $\ell^{2}(X)$ is called the \textit{uniform Roe algebra of $X$} and is denoted by $C^{\ast}_{u}(X)$.
\end{Def}

A representation of $\mathbb{C}[X]$ is a $*$-homomorphism $\pi: \mathbb{C}[X] \rightarrow \mb{B}(H)$, where $H$ is some Hilbert space. Each injective representation $\pi$ gives rise to a completion $C^{*}_{\pi}(X)\coloneqq \overline{\pi(\mb C[X])}\subset \mb B(H)$. In this context we think of $C^{*}_{u}(X)$ as the \textit{regular} completion. 

Using this observation, it is possible to show that a \textit{maximal} \Cs-norm makes sense and this leads to:

\begin{Def}
The maximal Roe algebra $C^{*}_{\max}(X)$ is the completion of $\mathbb{C}[X]$ in the norm %\footnote{We should probably just consider \textit{cyclic representations.}}
\begin{equation*}
\Vert T \Vert := \sup \lbrace \Vert \pi(T) \Vert \mid \pi \mbox{ a cyclic representation of } \mathbb{C}[X] \rbrace.
\end{equation*}
\end{Def}

\begin{Def}
Let $X$ be a coarse space with uniformly locally finite coarse structure $\mc E$, and let $E \in \mc E$ be an entourage. Then the \textit{E-Laplacian}, denoted by $\Delta^{E}$, is the element of $\mathbb{C}[X]$ with matrix entries defined by:
\begin{equation*}
\Delta^{E}_{x,y} = \begin{cases}
-1, & (x,y)\in (E\cup E^{-1})\setminus \diag(E),\\
\left| \lbrace z \in X \mid (x,z) \in (E\cup E^{-1})\setminus \diag(E)\right|, & x=y,\\
0 & \text{otherwise}. 
\end{cases}
\end{equation*}
Note that if $E \subset \diag(X)$ then $\Delta^{E}=0$. 
\end{Def}

\begin{Ex}\quad
\begin{enumerate}
\item If $X$ is a connected graph of bounded degree, then the set $E_{1}$, that is all pairs of points of distance $1$ (i.e the edges of the graph) generates the metric. In particular, $\Delta^{E_{1}}$ is the unnormalised graph Laplacian of $X$;
\item If $\Gamma$ is a finitely generated group, and then we can refine this above example to get the Laplacian:
\begin{equation*}
\Delta^{E_{1}}= 1-\sum_{s \in S}[s],
\end{equation*}
where $[s]$ is the formal element in the group ring $\mathbb{C}\Gamma$ given by $s \in S$, and $S$ (symmetrically) generates $\Gamma$ -- this \emph{group Laplacian} will be denoted by $\Delta_{\Gamma}$.
\end{enumerate}
\end{Ex}

This latter example connects with property (T) via a result of Valette \cite[Theorem 3.2]{MR773186}, which states that $\Gamma = \langle S\rangle$ has property (T) if and only if $0$ is isolated in the spectrum of the operator $\Delta_{\Gamma}$ in the maximal group $C^{*}$-algebra $C^{*}(\Gamma)$. 

Before moving onto the main result of this section, we point out that we can identify the algebraic Roe algebra, up to $*$-isomorphism, with the groupoid convolution algebra $C_{c}(G(X))$ \cite[Section 10.4]{MR2007488}, \cite[Appendix C]{Spakula-Willett}. In this way, groupoid reductions such as restricting to the boundary $\partial \beta X$ give rise to representations of $\mathbb{C}[X]$.

\begin{Def}
A representation of $\mathbb{C}[X]$ (or equivalently $C_{c}(G(X))$) is a \textit{boundary representation} whenever the ideal
\begin{equation*}
I_X = \lbrace T \in \mathbb{C}[X] \mid T_{x,y}\not =0 \mbox{ for only finitely many } x,y \in X \rbrace
\end{equation*}
is contained in the kernel.
\end{Def}
Note that in groupoid terms, $I_X$ is precisely the ideal $C_{c}(X\times X)$ in $C_{c}(G(X))$. Thus, a representation of $C_{c}(G(X)$ is a boundary representation if and only if it factors through $C_{c}(\partial G(X))$.

\begin{Def}
The \emph{boundary completion} $C^*_\partial(X)$ of $\mathbb{C}[X]$ is its separated completion in the seminorm
\[
\norm{T}_\partial \coloneqq \sup \lbrace \Vert \pi(T) \Vert \mid \pi \mbox{ a boundary representation of } \mathbb{C}[X] \rbrace
\]
\end{Def}

We can now state the relevant form of the definition of geometric property (T), using \cite[Proposition 5.2]{MR3246936}:

\begin{Def}
A space $X$ has \textit{geometric property (T)} (resp. \textit{geometric property (T) for boundary representations}) if there exists\footnote{This is equivalent to ``for every'' entourage, as the referenced proposition explains.} an entourage $E \in \mc E$ and a $c>0$ such that $\Spec_{\max}(\Delta^{E})$ (resp. $\Spec_{\partial}(\Delta^{E})$) is contained in $\lbrace 0 \rbrace \cup [c,\infty)$. Here $\Spec_\max$ denotes the spectrum in $C^*_\max(X)$ and $\Spec_\partial$ denotes the spectrum in $C^*_\partial(X)$.
\end{Def}

We note that the presence of the invariant measure $\mu$ on $\partial\beta X$ allows us to use the following well known $C^{*}$-algebraic fact:

\begin{Lemma}[{\cite[Section 7]{MR3246936}}]
Let $\Gamma\curvearrowright X$ be an action of $\Gamma$ on a compact Hausdorff space. Then $C^{*}_{\max}(\Gamma) \rightarrow C(X)\rtimes_{\max}\Gamma$ is injective if and only if $X$ has an invariant measure.\qed
\end{Lemma}
\begin{Cor}\label{cor:c*max-inj}
Let $\Gamma$ be a sofic group and $Z$ be a core of its sofic approximation. Then the natural map $C^{*}_{\max}(\Gamma) \rightarrow C(Z)\rtimes_{\max}\Gamma$ is injective.
\end{Cor}

\begin{Def}
We call any representation $\pi$ of $C_{c}(G(X))$ that factors through $C_{c}(\mc G|_{Z})$ \textit{sofic with respect to $Z$} or a \textit{$Z$-representation}. The \emph{sofic completion} $C^*_s(X)$ of $\mathbb{C}[X]$ is its completion in the norm
\[
\norm{T}_s \coloneqq \sup \lbrace \Vert \pi(T) \Vert \mid \pi \mbox{ a }  Z\mbox{-representation of } \mathbb{C}[X] \rbrace
\]
\end{Def}
Note that $C^*_s(X)\cong C^{*}_{\max}(\mc G|_{Z})$.

\begin{Def}
 $\mc X$ has \textit{geometric property (T) for $Z$-representations} if there exists $E \in \mc E$ and a $c>0$ such that $\Spec_s(\Delta^{E}) \subset \lbrace 0 \rbrace \sqcup [c,\infty)$, where $\Spec_s$ is the spectrum in $C^{*}_{s}(X)$.
\end{Def}

\begin{Thm}\label{Thm:T}
Let $\Gamma$ be a sofic group, $\mc X$ a sofic approximation and $X$ the corresponding space of graphs. Then $\Gamma$ has property (T) if and only if $X$ has geometric property (T) for $Z$-representations for any sofic core $Z\subset \partial\beta X$.
\end{Thm}
\begin{proof}
The proof is follows that of \cite[Theorem 7.1]{MR3246936}, making use of the fact that the operator $\Delta_{\Gamma}=\sum_{s\in S}1 - [s]\in \mathbb{C}\Gamma$ maps to the operator $\Delta^{Z} = \sum_{s \in S}1 - \sigma(s)$ in $C(Z)\rtimes_{\mathrm{alg}} \Gamma$, and thus it satisfies:
\begin{equation*}
\Spec_{\max}(\Delta_{\Gamma}) = \Spec_{\max}(\Delta^{Z}).
\end{equation*}
The result now follows from \cite[Theorem 3.2]{MR773186}, which shows that property (T) is equivalent to a spectral gap for $\Delta_{\Gamma}$.
\end{proof}

\begin{Cor}
If $\mc X$ has either geometric property (T) or geometric property (T) for boundary representations, then $\Gamma$ has property (T). \qed
\end{Cor}

%The following is an internal question that needs a little explanation: is the above corollary actually an equivalence? I would suspect that the answer should be yes, by some trickery involving the passage to a von Neumann setting. The following makes this a little more precise:

%If $G$ has property (T), then the groupoid $Z \rtimes G$ has property (T), and thus the groupoid $\mc G$ has property (T) as a \textit{measured} groupoid. This follows as $Z \rtimes G$ is an essential reduction of $\mc G$. Now, this gives us some (co)rigidity type statements for $\mc G$, from which we might be able to extract enough to get something like geometric (T) (or at least boundary geometric (T))?

\subsection{Locally embeddable into finite groups and some examples}\label{sect:LEF}

A group that is locally embeddable into finite groups has a $\eps=0$ sofic approximation $\mc X$, which we call an \textit{LEF approximation}. The set $Z$ in this case is the entire boundary $\partial\beta X$. From this we can observe that it is possible to prove the converse of some of the results from the previous section. This reproves essentially all of the results from \cite{Mimura-Sako} and \cite{MR3342685}. The arguments are straightforward after unpacking all of the definitions using groupoids.

\begin{Thm}
Let $\Gamma$ be LEF, let $\mc X$ be a LEF approximation and let $X$ be the space of graphs constructed as in section \ref{sect:ultralimits}. Then: 
\begin{enumerate}
\item $\Gamma$ is amenable if and only if $X$ has property A;
\item $\Gamma$ is property (T) if and only if $X$ has geometric property (T).
\end{enumerate}
\end{Thm}
\begin{proof}
It clearly suffices to prove the converses.

For i): as $\partial G(X)$ is topologically amenable, it has weak containment and a nuclear reduced groupoid $C^{*}$-algebra by \cite[Corollary 5.6.17]{brown-ozawa}. Additionally, the sequence
\begin{equation*}
0 \rightarrow \mc K (\ell^{2}(X)) \rightarrow C^{*}_{u}(X) \rightarrow C^{*}_{r}(\partial G(X)) \rightarrow 0,
\end{equation*}
is exact because of weak containment. It follows that $C^{*}_{u}(X)$ is nuclear, which is a well known characterisation of property A \cite{MR1905840}.

To show  ii), we immediately observe that every boundary representation is sofic, and hence boundary geometric property (T) follows. Moreover, the image of the group Laplacian $\Delta_{\Gamma}$ in representations of $G(X)$ given by convolution on the fibres of the source map retains spectral gap from property (T) by Corollary \ref{cor:c*max-inj}. This completes the proof.
\end{proof}

We remark that there are many interesting groups that are not residually finite, but are LEF -- chief amongst these are topological full groups of Cantor minimal systems, introduced by Giordano, Putman, and Skau \cite{MR1710743}, proved to be LEF by Grigorchuk and Medynets \cite{MR3241829}, amenable by Juschenko--Monod \cite{MR3071509} and have a simple commutator subgroup by Matui \cite{MR2205435}.

%I am aware these need to get a little longer and more precise but I have to think about how.
%\begin{Ex}
%(Full groups of Cantor minimal systems) Let $C$ be a Cantor set and let $T:C \rightarrow C$ be a homeomorphism. Then there is a group, denoted $[[T]]$, consisting of transformations of $C$ that locally look like powers of $T$. This group has been studied throughly and it was shown to be LEF \cite{MR3241829}, and that its commutator subgroup is amenable \cite{MR3071509}.

%Combining these results with the theorem above, we can conclude that the family of finite groups, which are set products of symmetric groups of certain sizes (or lamplighter-esque groups, see \cite{MR3314104}) form a family of groups with property A. This furnishes us with more examples of families of this type, adding to those of \cite{Mimura-Sako}.

%The main reason this example is interesting is precisely because $[[T]]$ is known to not be residually finite, so these techniques really do produce something completely new.
%\end{Ex}

\section{Coarse equivalence, quasi-isometry and uniform measure equivalence}\label{sect:coarseequiv}
In this section we prove that coarsely equivalent sofic approximations give rise to a uniform measure equivalence between groups, using Morita equivalence of groupoids as a tool. We first recall some definitions concerning the various notions of equivalence for groupoids that appear in the literature. 

\begin{Def}\label{def:linking}(A linking groupoid)
Let $G$ be a groupoid and let $T$ be a set with a map $f: T \rightarrow G^{(0)}$. Then the set
\begin{equation*}
G[T]:=\left\lbrace (t,t',g) \in (T\times T)\times G \mid g \in G^{f(t)}_{f(t')} \right\rbrace 
\end{equation*}
is a groupoid with the obvious operations. If $G$ is a locally compact Hausdorff groupoid, $T$ is a locally compact Hausdorff topological space and the map $f$ is continuous, then $G[T]$ is a locally compact Hausdorff topological groupoid. 
\end{Def}

For any sets $X,Y,T$ with maps $f:X \rightarrow T$, $g: Y \rightarrow T$ we denote the pullback by $X\times_{f,g} Y$, or $X \times_{T} Y$ if there is no ambiguity. 

\begin{Def}\label{def:action}(A groupoid action)
Let $G$ be a groupoid and let $M$ be a set. $M$ is a (right) $G$-space if there exists
\begin{itemize}
\item a map $p: M \rightarrow G^{(0)}$ (called the anchor map) and
\item a map $M \times_{p,r} G \rightarrow M$ denoted by $(z,g)\mapsto zg$ (called the action map)
\end{itemize}
with the following properties:
\begin{itemize}
\item $p(zg)= s(g)$ for all $(z,g) \in M \times_{p,r} G$;
\item $z(gh)=(zg)h$ whenever $p(z)=r(g)$ and $s(g)=r(h)$;
\item $zp(z) = z$ for every $z \in M$.
\end{itemize}
This allows us to define a natural ``crossed product'' groupoid $M\rtimes G$ with base space $M$, which consists of the elements $(z,z',g) \in (M\times M)\times G$ that satisfy $z=z' g$. Note that since $M\rtimes G \rightarrow M\times G$ given by $(z,z',g) \mapsto (z,g)$ is injective, we can also consider $M\rtimes G$ as a subset of $M \times G$, which we will do. We can also define a left $G$-space similarly using the source map instead of the range map: we denote the groupoid constructed from a left action by $G\ltimes M$.

Every groupoid $G$ naturally acts on its base space $G^{(0)}$ using $\id\colon G^{(0)}\to G^{(0)}$ as the anchor map and the multiplication as the action map. From the algebraic structure of the groupoid it easily follows that the orbit relation on $G^{(0)}$ defined by $x\sim y$ iff $x = y\cdot g$ for some $g\in G$ is an equivalence relation. The corresponding quotient is denoted by $G^{(0)}/G$. For a set $M$ with a $G$-action the quotient space by the action is defined through $M/G\coloneqq M/(M\rtimes G)$.
\end{Def}

So far we have mentioned nothing concerning the topological structure of the action and the crossed product space in case $G$ and $M$ have topologies. This can be adjusted by putting sufficient continuity and openness conditions on the maps above, which is discussed at length in \cite[Section 2]{MR2117427}. The main result of these considerations which we will need is the following:

\begin{Prop}\cite[Proposition 2.29]{MR2117427}\label{prop:morita_def}
Let $G_{1}$ and $G_{2}$ be two topological groupoids, let $s_{i},r_{i}$ be the open source and range maps of $G_{i}$. Then the following are equivalent:
\begin{enumerate}
\item there exists a set $T$ with $f_{i}:T \rightarrow G_{i}^{(0)}$ open surjective maps such that $G_{1}[T] \cong G_{2}[T]$;
\item there exists a space $M$ with two continuous maps $\rho: M \rightarrow G_{1}^{(0)}$, $\sigma: M \rightarrow G_{2}^{(0)}$ such that $\rho$ is the anchor map for a left action of $G_{1}$ on $M$, $\sigma$ is the anchor map of a right action of $G_{2}$ on $M$ such that these actions commute, are free and the action of $G_{2}$ is $\rho$-proper, the action of $G_{1}$ is $\sigma$-proper and:
\begin{equation*}
M/G_{2} \rightarrow G_{1}^{(0)} \mbox{ and } G_1 \backslash M \rightarrow G_{2}^{(0)} 
\end{equation*}
are homeomorphisms.
\end{enumerate}
Two topological groupoids that satisfy either of the two equivalent conditions above will be called  \textit{Morita equivalent}.
\end{Prop}

\begin{Rem}\label{Rem:bispace}
The main point to raise here is that the space $M$ in the proof of i) $\Rightarrow$ ii) is constructed as follows \cite[Proposition 2.29]{MR2117427}: take $M_{1}$ to be the space $G_{1}\times_{s,f_{1}}T$, and $M_{2}$ to be the space $T\times_{f_{2},r} G_{2}$. These are then combined over the $G_{i}[T]$-action on the right on $M_{1}$ and the left on $M_{2}$ to the space $M\coloneqq M_{1}\times_{G_{1}[T]}M_{2}$, which amounts of dividing the space $M_{1}\times_{T}M_{2}$ by the relation generated by $(z,z^{'}) \sim (zg, g^{-1}z)$, where $g \in G_{1}[T]$. The space $M$ then admits a bispace structure which implements ii).
\end{Rem}

\begin{Rem}
The notion of Morita equivalence can also be defined for measured groupoids in a similar manner, replacing topological conditions by measurable ones, and we will make use of it later. We refer the reader to \cite{landsman_quantized_2001} and references therein for discussion of definitions Morita equivalence for various categories of groupoids and operator algebras and connections between them.
\end{Rem}

\begin{Ex}
A coarse equivalence $f$ produces a ``coarse correspondence'', as in \cite{MR1905840}, between $G(X)$ and $G(Y)$. This is a groupoid $G(X\sqcup Y)$ constructed from a ``linking'' coarse structure, defined using the coarse structure $\mc E(f):=\mc E^{X}_{\mathrm{met}} \sqcup \mc E^{Y}_{\mathrm{met}} \sqcup \mc E^{XY} \sqcup \mc E^{YX}$, where the sets in $\mc E^{XY}$ are precisely those of the form $F\times f(F)$, similarly defining those in $\mc E^{YX}$ using the coarse inverse of $f$. This coarse structure is uniformly locally finite if $E^{X}_{\mathrm{met}}$ and $E^{X}_{\mathrm{met}}$ are \cite[Proposition 2.3]{MR1905840}.
\end{Ex}

This coarse correspondence allows us to construct a topological space $T=\beta X \sqcup \beta Y$ that implements a topological Morita equivalence between $G(X)$ and $G(Y)$ in the sense of Proposition \ref{prop:morita_def}. The proof of this is a part of the content of a remark from the beginning of Section 3.4 of \cite{MR1905840}.

\begin{Lemma}
If $X$ and $Y$ are coarsely equivalent by a pair of maps $f:X \rightarrow Y$ and $k:Y \rightarrow X$, then $G(X)[T] \cong G(Y)[T]$ for $T = \beta X$ and maps $p_{X}: T \rightarrow \beta X$ (resp. $p_{Y}:T \rightarrow \beta Y$) given by 
\begin{equation*}
p_{X}(\omega) = \begin{cases} \omega \mbox{ if } \omega \in \beta X \\ \overline{f}(\omega) \mbox{ if } \omega \in \beta Y \end{cases}.
\end{equation*}
and a similar definition for $p_{Y}$.\qed
\end{Lemma}

The space $M$ whose construction was outlined in Remark \ref{Rem:bispace} is a quotient of
\begin{equation}\label{eq:M}
M:=G(X)\times_{s,p_{X}}T\times_{p_{Y},r}G(Y)/\sim
\end{equation}
where $\sim$ implements the identification of points in $T$ who are joined by continuous extensions of the coarse maps $f: X\rightarrow Y$ and $k: Y \rightarrow X$. We also remark, that as the sets $X$ and $Y$ are invariant in their respective coarse groupoids, these bispaces restrict to bispaces over the boundary groupoids $\partial G(X)$ resp. $\partial G(Y)$.

\begin{Lemma}\label{lem:soficcoreCE}
Let $\Gamma$ and $\Lambda$ be sofic groups with $\mc X$, and $\mc Y$ sofic approximations of $\Gamma$ and $\Lambda$ respectively, and suppose that $f: X_{\mc X}\rightarrow X_{\mc Y}$ is a coarse equivalence of the associated spaces of graphs. Then the set $\widetilde{f}(Z_{\mc X})\cap Z_{\mc Y}$ has positive measure in $\partial\beta X_{\mc Y}$.
\end{Lemma}
\begin{proof}
By \cite[Lemma 1]{Khukhro-Valette} we can assume that $f(X_{i}) \subset Y_{i}$, and that $f|_{X_{i}}$ is a $(C,C)$-quasi-isometry (for some constant $C>0$). As $f$ is a coarse equivalence, there is a constant $n>0$ such that $X_{\mc Y} = N_{n}(f(X_{\mc X}))$, where $N_{n}$ is the $n$-neighbourhood of $f(X_{\mc X})$ in $X_{\mc Y}$. We also observe that $N_{m}(A) = N_{1}(N_{m-1}(A))$ for all subsets $A \subseteq X_{\mc Y}$ and all $m\in \mathbb{N}$. It follows by induction that, for all $i$:
\begin{equation*}
\vert N_{i}(f(X_{i})) \vert \leqslant \vert S_{\Lambda} \vert^{i} \vert f(X_{i}) \vert,
\end{equation*}
where $S_{\Lambda}$ is the finite generating set of $\Lambda$. This shows that $\widetilde{f}(Z_{\mc X})$ has measure at least $\frac{1}{\vert S_{\Lambda} \vert^{n}} \vert Y_{i} \vert$ in $\partial\beta Y$ as the image preserves unions and the measure $Z_{\mc X}^{c}$ is $0$. This completes the proof, since $Z_{\mc Y}$ has $\mu_{\mc Y}$-measure 1.
\end{proof}

In fact, we can say more using the observation that $\widetilde{f}(Z_{\mc X})\cap Z_{\mc Y} \not = \varnothing$: it allows us to construct a quasi-isometry using the transplanting technique of \cite[Proposition 3]{Khukhro-Valette}:

\begin{Prop}\label{prop:KV}
Let $\Gamma$ and $\Lambda$ be sofic groups, with sofic approximations $\mc X$ and $\mc Y$ respectively. If the spaces of graphs $X_{\mc X}$ and $X_{\mc Y}$ attached with $\mc X$ and $\mc Y$ are coarsely equivalent, then $\Gamma$ and $\Lambda$ are quasi-isometric.
\end{Prop}
\begin{proof}
The proof of this fact is precisely the proof of \cite[Proposition 3]{Khukhro-Valette}, except that instead of using convergence  of marked groups (i.e ultralimits of groups using the identity as base point), we use ultralimits along a base point sequence $(x_{i})_{i}$, such that $\eta = \lim_{\omega}x_{i}$ satisfies: $f(\eta) \in \widetilde{f}(Z_{\mc X})\cap Z_{\mc Y}$.
\end{proof}

Finally, we consider the analogous notion of measure equivalence, as was considered in \cite{Das-box} for box spaces of residually finite discrete groups. 

\begin{Def}[{\cite{MR1253544,MR2096453,Das-box}}]\label{def:ume}
Two groups $\Gamma$ and $\Lambda$ are \textit{measure equivalent} if there exists a essentially free Borel measure $(\Gamma,\Lambda)$-space $(X,\mu)$ such that there are finite volume fundamental domains $X_{\Gamma} \subset X \supset X_{\Lambda}$ for the actions. A measure equivalence is \textit{uniform} if additionally, for every $g\in \Gamma$ (resp. $h \in \Lambda$) there exists a finite subset $S_{g} \subset \Lambda$ (resp. $T_{h} \subset \Gamma)$ such that
\begin{equation*}
gX_{\Lambda} \subset X_{\Lambda}S_{g} \mbox{ and } X_{\Gamma}h\subset T_{h}X_{\Gamma}
\end{equation*}
\end{Def}

Our aim is to prove that if $\Gamma$ and $\Lambda$ are sofic groups with coarsely equivalent approximations, then $\Gamma$ and $\Lambda$ are uniformly measure equivalent. To accomplish this, we need to take the topological Morita equivalence $M$ of $G(X_{\mc X})$ and $G(X_{\mc Y})$ provided by a coarse equivalence $f: X_{\mc X} \rightarrow X_{\mc Y}$, and turn it into a Morita equivalence between measured groupoids. To do this, we have to analyse the correspondence between invariant measures and measures on a quotient by a free and proper action for \'etale groupoids:
\begin{Prop}\label{prop:morita}
Let $G$ and $H$ be \'etale groupoids and let $X$ be a free and proper $G$-$H$-space. Then there is a one-to-one correspondence between $G$-invariant Radon measures $\rho$ on $X$ and Radon measures $\mu$ on $G\backslash X$. Moreover, this correspondence is additive and $H$-equivariant.
 % and $\mu_H$ on $G^{(0)}$ and $H^{(0)}$ respectively, and let $M$ be a Morita equivalence. Then there is a $G$-$H$-bi-invariant measure on $M$.
\end{Prop}

 % construct a bi-invariant measure on the space $M$ in the situation that there are invariant measures on the boundaries $\partial\beta X_{\mc X}$ and $\partial \beta X_{\mc Y}$. We do this in general for \'etale groupoids:

\begin{proof}
Each $G$-invariant Radon measure $\rho$ on $X$ defines a Radon measure $\mu=\leftsub{G}{\rho}$ on $G\backslash X$ using the pushforward of $\rho$ over a subset $U\subset X$ such that the quotient map is one-to-one on $U$. This construction is $H$-equivariant as the $H$-action commutes with the $G$-action on $X$.

To go back, we use the construction from \cite[Section 3]{MR2966043}: let $X$ be a free and proper left $G$-space. Then $G\backslash X$ is a locally compact Hausdorff space, and for each $x \in X$, the map $\gamma\mapsto \gamma\cdot x$ is a homeomorphism of $G\cdot r^{-1}(x)$ onto the orbit $G\cdot x$. We define a Radon measure $\rho^{G\cdot x}$ on $X$ with support $G\cdot x$ by
\[
\rho^{G\cdot x}(f):=\int_G f(\gamma^{-1}(x))d\lambda^{r(x)}(\gamma)
\]

Our definition is independent of our choice of $x$ in its orbit by left-invariance of the
Haar system $\lambda$. Additionally, the map
\[
[x]\mapsto \rho^{[x]}(f)
\]
is continuous on $G\backslash X$. Given a finite Radon measure $\mu$ on $G\backslash X$, we define a Radon
measure $\rho_\mu$ on $X$ by
\[
\rho_\mu(f) = \int_{G\backslash X} \int_X f(y) d\rho^{[x]}(y) d\mu([x])
\]

% Applying this construction to $M$ and noticing that the Haar measure on $r^{-1}(x)$ is the counting measure as well as using the isomorphism $G\backslash M \cong H^{(0)}$, we obtain
% \[
% \rho(f) = \int_{H^{(0)}} \int_X f(y) d\rho^{h}(y) d\mu_H(h)
% \]
The measure $\rho$ is $G$-invariant by construction, as $\rho^{[x]}$ is invariant and supported on a $G$-orbit. On the other hand, as the actions of $G$ and $H$ on $X$ commute and because the measures $\rho^{[x]}$ are defined by integrating over the orbit, they are $H$-equivariant: for all $\chi\in H$ we have $\rho^{[x]\cdot\chi} = \chi_*\rho^{[x]}$.

It's routine to check that these constructions are additive, inverse to each other and therefore define a one-to-one correspondence as claimed.
\end{proof}

In the situation of the above proposition we say that $\mu$ is the quotient measure corresponding to $\rho$ and write $\mu = \leftsub{G}{\underline{\rho}}$ and that $\rho$ is the measure induced by $\mu$ through the action of $G$ and write $\rho = \leftsup{G}{\overline{\mu}}$; we use corresponding notations for right actions.

\begin{Cor}\label{cor:measurable-morita}
Let $G$ and $H$ be \'etale groupoids with invariant measures $\mu$ and $\eta$ on $G^{(0)}$ and $H^{(0)}$ respectively and let $M$ be a Morita equivalence between them. If $\underline{\leftsup{G}{\overline{\mu}}}_H$ on $H^{(0)}$ is absolutely continuous with respect to $\eta$ and $\leftsub{G}{\underline{\overline{\eta}^H}}$ on $G^{(0)}$  is absolutely continuous with respect to $\mu$, then $(G,\nu_\mu)$ and $(H,\nu_\eta)$ are Morita equivalent as measured groupoids in the sense of \cite{landsman_quantized_2001}.
\end{Cor}
\begin{proof}
The absolute continuity assumptions imply that the measure
\[
\rho\coloneqq \leftsup{G}{\overline{\mu}} + \overline{\eta}^H
\]
descends to measures $\leftsub{G}{\underline{\rho}}$ and ${\underline{\rho}}_H$ which are equivalent to $\mu$ resp. $\eta$. Thus, $(M,\rho)$ is a Morita equivalence between the measured groupoids $(G,\nu_\mu)$ and $(H,\nu_\eta)$.
\end{proof}
% This enables us to transform a topological Morita equivalence between two \'etale groupoids $G$ and $H$ into a measurable Morita equivalence between the groupoids $(G,\nu_{\mu_{G}})$ and $(H,\nu_{\mu_{H}})$.

Using this we can prove:

\begin{Thm}\label{thm:coarse-equiv}
Let $\Gamma$ and $\Lambda$ be sofic groups with approximations $\mc X$ and $\mc Y$ respectively. If the associated spaces of graphs $X_{\mc X}$ and $X_{\mc Y}$ are coarsely equivalent, then the groups $\Gamma$ and $\Lambda$ are uniformly measure equivalent.
\end{Thm}
\begin{proof}
In order to appeal to Corollary \ref{cor:measurable-morita}, we have to show that the limits $\mu$ and $\eta$ of counting measures on the base spaces of $\mc G_{\Gamma}$ and $\mc G_{\Lambda}$ satisfy the absolute continuity assumption. To check this, recall the construction of the space $M$ following \eqref{eq:M}:
\[
M\coloneqq \partial G(X)\times_{s,p_{X}}T\times_{p_{Y},r}\partial G(Y)/\sim
\]
where $\sim$ implements the identification of points in $T = \partial\beta X \sqcup \partial\beta Y$ who are joined by continuous extensions of the coarse maps $f: X\rightarrow Y$ and $k: Y \rightarrow X$. It follows that the measure $\underline{\leftsup{G}{\overline{\mu}}}_H$ is equal to the pushforward $\overline f_\ast \mu$ of the measure $\mu$ under the coarse equivalence map $f$, and similarly, $\leftsub{G}{\underline{\overline{\eta}^H}}$ is equal to the pushforward $k_\ast \eta$ under the coarse inverse. As coarse maps have uniformly finite fibres, the absolute continuity follows. Thus, Corollary \ref{cor:measurable-morita} yields a measurable Morita equivalence $(M,\rho)$ between $(\mc G_\Gamma,\nu_\mu)$ and $(\mc G_\Lambda,\nu_\eta)$. 

% Notice that  if $\mc G_{\Gamma}$ and $\mc G_{\Lambda}$ are topologically Morita equivalent, then there is a measurable $(\Gamma,\Lambda)$-bispace $M$ using Proposition \ref{prop:morita}, as Theorem \ref{thm:clevertrick2} shows us that  $\mc G_{\Gamma}$ and $\mc G_{\Lambda}$ are equal to $(\partial \beta X_{\mc X}, \mu_{\Gamma})\rtimes \Gamma$ and $(\partial\beta X_{\mc Y},\mu_{\Lambda})\rtimes \Lambda$ up to measurable isomorphism. Whence the space $(M,\rho)$ is a $(\Gamma,\Lambda)$-Borel measure space on which these groupoids act measurably, with compact (and finite volume) fundamental domains. This implies that $\Gamma$ and $\Lambda$ are measure equivalent. 

To show that $\Gamma$ and $\Lambda$ are uniformly measure equivalent, we fix fundamental domains $X_{\Gamma}, X_\Lambda \subset M$ with compact closures for the $\mc G_{\Gamma}$ and $\mc G_\Lambda$-actions respectively and let $\lbrace U_{g} \rbrace_{g\in\Gamma}$ and $\lbrace U_{h} \rbrace_{h\in\Lambda}$ be covers of $\mc G_{\Gamma}$ and $\mc G_\Lambda$ by compact open slices, each of which restricts to a slice of the form $[Z_{\mc X},g]$ on $\mc G|_{Z_{\mc X}}\cong Z_{\mc X} \rtimes \Gamma$ and $[Z_{\mc Y},h]$ on $\mc G|_{Z_{\mc Y}}\cong Z_{\mc Y} \rtimes \Lambda$.
 
Fix $h\in\Lambda$. The set $\lbrace U_{g}X_{\Lambda} \rbrace_{g \in \Gamma}$ is an open cover of $M$, thus in particular it covers $X_{\Gamma}U_{h}$, which is a compact subset of $M$ as the right $\mc G_{\Lambda}$ action is proper and $U_{h}$ is a compact open slice of $\mc G_{\Lambda}$. Now, compactness of $X_{\Gamma}U_{h}$ allows us to extract a finite subcover $\lbrace U_{g}X_{\Gamma} \rbrace_{g\in T_{h}}$ for some finite set $T_{h}\subset \Gamma$.

To finish the proof, we remark that the the almost everywhere isomorphisms $\mc G_\Gamma \cong \partial\beta X\rtimes \Gamma$ and $\mc G_\Gamma \cong \partial\beta X\rtimes \Lambda$ constructed in Theorem \ref{thm:ae-iso} give rise to actions of $\Gamma$ and $\Lambda$ on $M$ with (measurable) fundamental domains $X_\Gamma$ and $X_\Lambda$ such that $g X_\Gamma$ and $X_\Gamma h$ coincide with $U_gX_\Gamma$ and $X_\Gamma U_h$ up to null sets. Thus, the set $T_h$ satisfies the condition in the Definition \ref{def:ume}, and symmetrisation of the  argument for the $\mc G_{\Lambda}$-action provides for evey $g\in G$ a finite set $S_{g}$ with the necessary properties. This finishes the proof.
\end{proof}

Appealing to \cite[Theorem 6.3]{MR1953191}, we now obtain: 

\begin{Cor}
If $\Gamma$ and $\Lambda$ are finitely generated sofic groups with coarsely equivalent sofic approximations, then their $\ell^{2}$-Betti numbers are proportional.
\end{Cor}

This Corollary has immediate applications to distinguishing families of finite graphs up to coarse equivalence. In particular, it allows us to see that box spaces of products of free groups with different number of factors are not coarsely equivalent \cite[Corollaire 0.3]{MR1953191}) as they have $\ell^{2}$-Betti numbers that are not proportional -- we remark that this is considered directly in the work of Das \cite{Das-box}, and we draw attention to it again due to recent interest in this question \cite{Khukhro-Valette,fullboxspace}.
  
\subsection{Remarks about bilipschitz equivalence}
Let $\Gamma$ and $\Lambda$ be sofic groups with approximations $\mc X$ and $\mc Y$ respectively. If $X_{\mc X}$ and $X_{\mc Y}$ are bilipschitz equivalent via a map $f$, then they are certainly coarsely equivalent and so the results of the previous section apply. However, as in the remark that precedes \cite[Definition 2.1.4]{MR2096453}, we can say quite a bit more concerning the relationship between $\Gamma$ and $\Lambda$ in this instance.

Notably, the following basic observations can be used to simplify and improve on the results from Section \ref{sect:coarseequiv}:

\begin{enumerate}
\item Bilipschitz equivalences are bijections, so the pushforward $f_{*}\mu_{\mc X}$ agrees with $\mu_{\mc Y}$. This means that Lemma \ref{lem:soficcoreCE} is a triviality, as $\tilde{f}(Z_{\mc X})$ is $\mu_{\mc Y}$-measure 1. We also remark that any bijection from $X_{\mc X}$ to $X_{\mc Y}$ will also give a homeomorphism between $\partial \beta X_{\mc X}$ and $\partial \beta X_{\mc Y}$;
\item Let $\mu$ and $\eta$ be measures on $G^{(0)}$ and $H^{(0)}$ respectively (as in Corollary \ref{cor:measurable-morita}). Then applying i), but this time in the construction of the bimodule measure $\rho$ induced from $\mu$, we see that actually $\underline{\leftsup{G}{\overline{\mu}}}_H=\eta$. As a consequence,
\begin{enumerate}
\item we do not need to use the sum of $\rho:=\leftsup{G}{\overline{\mu}} + \overline{\eta}^H$ in the proof of Corollary \ref{cor:measurable-morita};
\item there is a common fundamental domain in a topological and measurable sense (as a consequence of the homeomorphism between $\partial \beta X_{\mc X}$ and $\partial \beta X_{\mc Y}$).
\end{enumerate}
\end{enumerate}

Additionally, one can improve Proposition \ref{prop:KV}.

\begin{Prop}
Let $\Gamma$ and $\Lambda$ be sofic groups, with sofic approximations $\mc X$ and $\mc Y$ respectively. If the spaces of graphs $X_{\mc X}$ and $X_{\mc Y}$ attached with $\mc X$ and $\mc Y$ are bilipschitz equivalent, then $\Gamma$ and $\Lambda$ are bilipschitz equivalent. \qed
\end{Prop}

This has additional consequences due to results by Medynets--Thom--Sauer \cite[Theorem 3.2]{medynets-thom-sauer}:

\begin{Cor}
Let $\Gamma$ and $\Lambda$ be sofic groups, with sofic approximations $\mc X$ and $\mc Y$ respectively. If the spaces of graphs $X_{\mc X}$ and $X_{\mc Y}$ attached with $\mc X$ and $\mc Y$ are bilipschitz equivalent, then there exists minimal, continuous orbit equivalent actions of $\Gamma$ and $\Lambda$ on some Cantor set $C$.\qed
\end{Cor}

\section{A standardisation of the base space}\label{sect:standard}
This section is dedicated to the proof of the following theorem:
\begin{Thm}\label{thm:standard}
Let $\Gamma$ be a sofic group, $\mc X$ be a sofic approximation of $\Gamma$ and $X$ the associated total space of the family of graphs attached to $\mc X$. Then there exists a second countable \'etale, locally compact, Hausdorff topological groupoid $\mc G$ with following properties:
\begin{enumerate}
\item the base space $\mc G^{(0)}\eqqcolon \widehat X$ is a compactification of $X$ (in particular, it's a quotient of $\beta X$ through a quotient map $p\colon \beta X \to \widehat X$),
\item $p(Z)\subset \partial \widehat X$ is invariant and satisfies $\mc G|_{p(Z)} \cong p(Z)\rtimes \Gamma$. As a consequence, we have an almost everywhere isomorphism
\[
(\mc G|_{\partial \widehat X}, \nu_{p_{*}\mu}) \rightarrow (\widehat X,p_{*}\mu)\rtimes \Gamma.
\]
\end{enumerate}
\end{Thm}
Morally, this means that although the space $(\partial\beta X,\mu)$ is not a standard probability space, we can use $X_{\mc A}$ to make arguments \textit{as if} we were actually in $\partial\beta X$, whilst actually working in a standard Borel probability space.

\begin{Ex}
Let $\Gamma$ be a residually finite, finitely generated discrete group, let $\mc X$ be a sofic approximation made up of finite quotients of $\Gamma$ and let $X$ be the space of graphs associated to $\mc X$. Then by considering the Boolean algebra $B$ generated by $\Cofin(X) \cup \Fin(X) \cup \lbrace \Sh(eN_{i}) \rbrace_{i}$, where
\begin{equation*}
\Sh(eN_{i}) \coloneqq \bigcup_{j\geqslant i}\lbrace x \in X_{j} \mid \pi_{i,j}(x)=e_{i}N_{i} \rbrace 
\end{equation*}
is the \textit{shadow} of $e_{i}$ in $X$, we obtain a second countable, locally compact, Hausdorff \'etale groupoid $\mc G_{B}$, which is homeomorphic to $X_{B} \rtimes \Gamma$, and $X_{B} \cong \widehat{\Gamma}_{\mc X}$ is the profinite completion associated with the family of finite quotients $\mc X$. This dynamical system was introduced in \cite{MR2966663}, where it was shown to be minimal (and in this case, as subgroups in question are normal, it's also free). A similar construction using the shadows of the identity would give us the boundary $\partial T$ as defined in \cite{MR2966663} when the chain is \textit{Farber}. This example shows that one can choose the appropriate Boolean algebra depending on the goals in question. 
\end{Ex}

The ideas used in the proof stem from the work of Skandalis--Tu--Yu \cite{MR1905840}, where one pushes the failure of second countability of $G(X)$ purely into the unit space: this allows one to make use of the groupoid equivariant KK-theory of Pierre-Yves Le Gall \cite{MR1855247} to describe the coarse Baum--Connes conjecture attached to $X$.

\begin{proof}[Proof of Theorem \ref{thm:standard}]
 We first recall the outcome of \cite[Lemma 3.3]{MR1905840}, which states that any countable generating set $\mc A$ of the metric coarse structure on a space $X$ gives rise to a second countable, \'etale, locally compact Hausdorff groupoid $G_{\mc A}$, such that the coarse groupoid $G(X)$ is homeomorphic to the transformation groupoid $\beta X \rtimes G_{\mc A}$\footnote{Skandalis-Tu-Yu give a ``by hand'' proof of this result: a slightly more modern approach to it would be to make use of the fact that a pseudogroup in this context gives us a inverse monoid, and then construct from that, using well known techniques of \cite{MR2419901}, a groupoid with all the desired properties.}. We construct $\mc A$ in what follows.

% We however can get further traction from these ideas by first choosing a good generating set for the metric coarse structure, and then observing that the natural (surjective, open, continuous) anchoring map $p:\beta X \rightarrow X_{\mc A}$, where $X_{\mc A}$ is the unit space of $G_{\mc A}$, allows us to push forward the measure $\mu$ defined in Section \ref{sect:ultralimits} to $X_{\mc A}$.

In light of Section \ref{sect:sofic-coarse-bd-groupoid}, we can construct generators using those given by the labelling, i.~e. by considering the entourages $E_{P}$, where $P$ is a word in the free group on the alphabet $S$. These clearly generate the metric for the space $X$ (as a total space of the family $\mc X$). This family doesn't give us a good unit space however, as each of the elements we are using here are bijections (thus the base space $X_{\mc A}$ of $G_{\mc A}$ would end up being a point). 

To remedy this, we consider the set $\mc B$ of all countable Boolean subalgebras of $\textbf{2}^{X}$ that contain the set $Y$ and some infinite set that is not cofinite. Note that if the approximation $\mc X$ is a LEF approximation (i.e $\eps=0$ for $i$), then $Y=X$ and subsequently, this is all countable Boolean subalgebras with at least one infinite, not cofinite set.  

Fix $B \in \mc B$. By taking the inverse semigroup generated by $B$ and the transformations $\tau(w)$ for $w \in F_{S}$, we get a countable pseudogroup.
%(the process of closing to get an inverse monoid here ends up producing all the necessary restrictions, and also equivariantises the Boolean algebra in question).
Let $\mc A$ be this set of partial transformations $\tau(g)|_{A}$, where $A \in \mc B$, extended continuously to $\beta X$. Applying \cite[Lemma 3.3]{MR1905840}, we obtain a second countable \'etale groupoid $G_{\mc A}$. Its base space $X_{\mc A}$ is a quotient of $\beta X$ and we denote the quotient map by $p\colon \beta X\to X_{\mc A}$. Pushing forward the measure $\mu$ along the map $p\colon \beta X\to X_{\mc A}$, and using the Urysohn metrization theorem, we obtain that $(X_{\mc A},p_{*}\mu)$ has the structure of a standard Borel probability measure space. Since $\mu$ is supported inside the boundary $\partial\beta X$, it follows that $p_{*}\mu(\partial X_{\mc A})=1$, and we again obtain a boundary type groupoid $\mc G_{\mc A}:=G_{\mc A}|_{\partial X_{\mc A}}$, and as the set $Y$ from Section \ref{sect:sofic-coarse-bd-groupoid} belongs to the Boolean algebra generating $\mc G_{\mc A}$, we can see that the $F_{S}$-action factors through $\Gamma$ up to null sets. This allows us to run the arguments of Section \ref{sect:sofic-coarse-bd-groupoid} again to obtain an almost everywhere isomorphism of groupoids (through a $\mu$-inessential reduction). This finishes the proof.
\end{proof}

%\subsection{Cost}

%\begin{Def}
%Cost of an equivalence relation, combinatorial cost of a family, some connective tissue between these
%\end{Def}

%\begin{Thm}
%Let $\Gamma$ be a sofic group, let $\mc X$ be a sofic approximation and let $\mc G$ be the boundary groupoid associated to the space of graphs $X$ of $\mc X$. Then the combinatorial cost $cc(\mc X)$ is equal to the cost of the action $\Gamma \curvearrowright \partial \beta X$. 
%\end{Thm}

\section{Concluding remarks and further questions}
We finish the paper with a few questions and comments on the surrounding literature, concerning primarily the interactions between the geometric and probabilistic points of view on sofic groups and graphs. Throughout, let $\Gamma$ be a sofic group, $\mc X$ a sofic approximation and $\mc G|_{Z}$ be the sofic coarse groupoid restricted to the sofic core.

The statement of our main result immediately suggests a question about the converse:
\begin{Qu}
To which extent do the converse statements to the one of Theorem \ref{thm:intro} hold?
\end{Qu}

Because soficity gives only a measure-theoretic control of actions on the sofic boundary, we do not expect the converse to hold true in full generality. On the other hand, as amenability, a-T-menability and property (T) of discrete groups are visible at the level of measure-preserving actions, it is natural to expect that they will be visible at the sofic boundary; it is natural to expect some form of probabilistic manifestation of coarse-geometric properties there.

\begin{Def}[{\cite{MR2455943,MR2372897}}]
A family of finite graphs $\mc Y = \lbrace Y_{i} \rbrace_{i}$ of bounded degree is a \textit{hyperfinite family} if for every $\eps > 0$ and for each $i \in \mathbb{N}$ there exists a decomposition of $Y_{i}$ into $K_{\epsilon,i}$ finite sets $U_{i,j}$ such that
\begin{enumerate}
\item each $U_{i,j}$ is uniformly bounded;
\item the size of each set $E(U_{i,j},U_{i,j^{'}})$ is at most $\eps \vert Y_{i} \vert$ whenever $j \not = j^{'}$, where $E(U_{i,j},U_{i,j^{'}})$ is the set of edges between $U_{i,j}$ and $U_{i,j^{'}}$.
\end{enumerate}
\end{Def}

A combination of Theorem \ref{Thm:amenable} with \cite[Theorem 1.1]{MR2372897} shows that property A for a sofic approxiation implies hyperfiniteness of that approximation. 

\begin{Qu}
Does hyperfiniteness of a sofic approximation imply property A for that approximation?
\end{Qu}

Here the fact that we ask this for a sofic approximation is important, as the implication does not hold for a general Benjamini--Schramm convergent sequence of graphs\footnote{This was communicated by Gabor Elek, in a personal communication.}.

One approach to this question would be to use a property equivalent to property A called the \textit{metric sparsification property} \cite{MR2419930}, which was shown to be equivalent in \cite{MR3200343} to a graph family being weighted hyperfinite (as defined in by Elek and Tim\'{a}r in \cite{elek-timar}). However, there is a subtlety here -- the measure on the groupoid $\mc G$ only deals with counting measures on the graphs, whereas the weighted notion of hyperfiniteness from \cite{elek-timar} is dealing with limits of arbitrary measures on the graph family $\mc X$. 

Another recent development in \cite{sofic_T} classified measurably those approximations coming from groups with property (T). The natural analogue of hyperfiniteness in this setting is the following:

\begin{Thm}[{\cite[Theorem 1]{sofic_T}}]
Let $\Gamma$ be a property (T) group and let $\mc X=\lbrace X_{i} \rbrace_{i}$ be a family of bounded degree graphs that Benjamini-Schramm converge to the Cayley graph of $\Gamma$. Then there is a $c>0$ and a family of regular graphs $\mc Y = \lbrace Y_{i} \rbrace_{i}$ such that:
\begin{enumerate}
\item $V(X_{i})=V(Y_{i})$ for every $i$;
\item $\lim_{i \rightarrow \infty}\frac{\vert E(X_{i})\bigtriangleup E(Y_{i}) \vert}{\vert V(X_{i}) \vert} = 0$;
\item Each $Y_{i}$ is a vertex disjoint union of $c$-expanders.
\end{enumerate}
\end{Thm}

In other words, the graphs $Y_i$ are obtained by ``rewiring'' $X_i$ in an asymptotically negligible manner.

\begin{Qu}
Can a sofic approximation of a property (T) group be ``asymptotically rewired'' to have some form of geometric property (T)?
\end{Qu}

We remark that a combination Theorem \ref{Thm:T} implies that geometric (T), boundary geometric (T) or geometric (T) for sofic representations imply the conclusion of \cite[Theorem 1]{sofic_T}, so the above is asking about a strengthening of the latter.

As the results of \cite{sofic_T} are statements about the ergodic decomposition of the measure, and these specific questions motivate the following:

\begin{Qu}[Ergodic decomposition]\label{qu:ergodic}
What properties do the subgroupoids of $\mc G|_{Z}$ that correspond to the ergodic components have?
\end{Qu}
Notice that this question connects very nicely to older results, notably \cite{MR2455943} and \cite{MR2372897}. 

On a related note, there are many measurable notions from the literature, such as cost \cite{MR2455943}, entropy and mean dimension \cite{MR3130315} that all apply to \textit{measured} groupoids -- the topological groupoid defined in Section \ref{sect:sofic-coarse-bd-groupoid} can also be considered in this setting, and after passing through the standardisation process of Section \ref{sect:standard} we obtain groupoids that allow for these notions to be applied. This mirrors the work of Carderi \cite{carderi-ultraproducts}, as remarked earlier.

Our standardisation process produces topological groupoids, but is far from giving a unique space -- the difference being that we use countable Boolean subalgebras of $\textbf{2}^{\mc X}$, as opposed to countable Borel $\sigma$-algebras -- and these each give potentially give rise to very different metrisable dynamical systems. On the other hand, properties such as amenability and property (T) will pass to these systems without any loss. This naturally leads to the following question:

\begin{Qu}
What is the interaction between coarse properties of $\mc X$ and the measurable properties of its various standardisations? More concretely, can we show that for these systems, the invariants such as entropy (or mean dimension) do not depend on the choice of countable Boolean subalgebra? Is there an ``clopen" analogue of the main results in \cite{carderi-ultraproducts}?
\end{Qu}

% These questions are related doing an analogue of ``measurable'' coarse geometry: 

% \begin{Qu}
% What can we say about the von Neumann algebra of the groupoid $\mc G$ in general? How much of the theory can be profitably considered in this case?
% \end{Qu}

Finally, we end this section with a remark about a specific sofic group that itself does not have property A.

\begin{Ex}[Non-exact groups that are sofic]
It is known, by a construction proposed in \cite{arzhantseva-osajda} and completed in \cite{Osajda-nonexact} that there are groups that are a-T-menable, but do not have property A. A natural observation is that any such $\Gamma$ is direct limit of hyperbolic, CAT(0)-cubulable groups $\Gamma_{m}$ -- and as hyperbolic CAT(0)-cubical groups are residually finite \cite{MR3104553}, $\Gamma$ will be LEF (see \cite{MR3074498} for a proof of this, in the more general sofic setting).

In this situation, any LEF sequence will mostly likely be asymptotically coarsely embeddable, but it will not satisfy a notion of ``asymptotic property A'' that will be introduced in \cite{pillon}, which is some form of groupoid exactness that appears to fail in the general setting -- this is related to doing coarse geometry on groupoids with metrizable range fibres as in \cite{groupoidroealgebra} or \cite{exactgroupoids}.
\end{Ex}

\begin{Qu}
What can we say concerning the asymptotic geometry of the sofic approximations of the above monster groups? Can we use embeddings of sofic approximations to construct new exotic monster groups with strange properties?
\end{Qu}

\section*{Acknowledgements}
This research was initiated at the Erwin Schrödinger Institute programme ``Measured group theory", February 2016. Correspondingly, the authors would like to thank the ESI for the nice working atmosphere support during the workshop. During this period, both the authors were also partially supported by the ERC grant “ANALYTIC” no. 259527 of Goulnara Arzhantseva. The authors would additionally like to thank Rufus Willett for ideas that lead to some of the statements concerning ultralimits and amenability (particularly Proposition \ref{prop:withRufus}), Erik Guentner for an excellent discussion about ultralimits that reinvigorated the work, Andreas Thom and Ben Hayes for stimulating discussions and comments on an early version of this paper. We also thank the anonymous referee for valuable comments which helped us to improve the exposition.

\bibliographystyle{alpha}
%\bibliography{Remote}
\bibliography{sofic}

\end{document}